\documentclass{article}

\usepackage[round]{natbib}
\usepackage{fancyvrb}
\usepackage{amssymb}

\title{MATHEMATICS AND LOGIC AS INFORMATION\\
COMPRESSION BY MULTIPLE ALIGNMENT,\\
UNIFICATION AND SEARCH}

\author{J Gerard Wolff
\\
\\
{\small \it CognitionResearch.org.uk.}
\\
{\small \it E-mail: gerry@cognitionresearch.org.uk.}
}

\normalsize

\begin{document}

\maketitle

\begin{abstract}
This article introduces the conjecture that {\it mathematics, 
logic and related disciplines may usefully be understood as information compression by `multiple alignment', `unification' and `search'}.

As a preparation for the two main sections of the article, concepts of Hartley-Shannon information theory and Algorithmic Information Theory are briefly reviewed together with a summary of the essentials of the simpler `standard' methods for information compression. Related areas of research are briefly described: philosophical connections, information compression (IC) in brains and nervous systems, and IC in relation to inductive inference, Minimum Length Encoding and probabilistic reasoning. Then the concepts of information compression by `multiple alignment', `unification' and `search' (ICMAUS) are outlined together with a brief description of the SP61 computer model that is a partial realisation of the ICMAUS framework.

\sloppy{The first of the two main sections describes how many of the commonly-used forms and structures in mathematics, logic and related disciplines (such as theoretical linguistics and computer programming) may be seen as devices for IC. In some cases, 
these forms and structures may be interpreted in terms of the ICMAUS framework.}

The second main section describes a selection of examples where processes of calculation and inference in mathematics, logic and related disciplines may be understood as IC. In many cases, these examples may be understood more specifically in terms of the ICMAUS concepts and illustrated with output from the SP61 model.

Associated issues are briefly discussed.
\\
\noindent {\it Abbreviations used in the article:} IC---information compression; 
ICMAUS---information compression by multiple alignment, unification and search; 
MA---multiple alignment; ML---mathematics and logic; 
MLE---Minimum Length Encoding; OOD---object-oriented design; 
PCS---Post Canonical System; UTM---Universal Turing Machine.

\end{abstract}

\section{\bf Introduction}

This article introduces some conjectural ideas about the foundations of mathematics, logic and related disciplines such as theoretical linguistics and computer programming. In this context, `mathematics' means the the subject itself rather than mathematical method. For the sake of brevity, the expression `mathematics, logic and related disciplines' will be referred to as `ML'\footnote{Not to be confused with the programming language of the same name.}.

In broad terms, the article examines how {\em information theory} may illuminate the nature of ML, an approach to the subject that seems to have received little or no attention in the past. Here, `information theory' encompasses Shannon's \citeyearpar{shannon_weaver_1949} information theory and the more recently developed `algorithmic information theory' \citep[see, for example,][]{li_vitanyi_1997}.

In brief, the main conjecture to be discussed is the proposition that:

\begin{quotation}

\noindent {\it Mathematics, logic and related disciplines may usefully be understood as information compression by `multiple alignment', `unification' and `search'.}

\end{quotation}

\noindent The concepts of {\it information compression} (IC), {\it multiple alignment} (MA), {\it unification} and {\it search} (collectively abbreviated as ICMAUS) will be described.

As will be seen, the conjecture will be discussed in relation to the `static' representation of concepts in ML, and also in relation to the `dynamics' of ML inferences: `calculation', `proof', `deduction' etc.

For many readers, it may seem strange to suppose that there should be any 
significant connection between ML and IC. After all, IC is the function performed by those humble but useful computer utilities like PkZip and WinZip that can be used to reduce the sizes of files but do not otherwise seem to have any particular significance. In this article, I hope to show that there are many other facets to IC than what is apparent in compression utilities and that there are deep intellectual roots to the ideas to be described.

There are no formal proofs in this article and indeed formal proofs in this context 
are probably inappropriate since the argumentation depends on one's sense of analogy and 
what is or has potential to be intellectually fruitful. The treatment of the 
proposals is very far from being exhaustive: the article merely introduces 
some ideas in support of the proposals and presents a selection of examples. Whether or how the ideas may be applied to the many areas of ML not discussed in this article are matters that readers may like to consider.

\subsection{`COMPUTING', ML AND IC}\label{computing_ML_IC}

In another article \citep{wolff_1999_prob}, I have suggested that the intuitive concept of 
`computing' may be understood as ICMAUS. If the arguments in that article are accepted, 
and if `computing' is seen to include ML, then the present article might be thought to be 
redundant.

However, the purpose of the previous article is to show that the workings of a 
Universal Turing Machine (UTM) and the workings of a Post Canonical System (PCS) may 
be understood in terms of ICMAUS and to discuss how this interpretation suggests how the 
UTM and PCS models may, with profit, be augmented so that the practical benefits of the 
augmentation outweigh the costs. By contrast, the purpose of this article is to show how 
specific aspects of ML other than UTMs or PCSs may be understood in terms IC and, more 
specifically in some cases, in terms of ICMAUS. The two articles may be seen to be 
complementary rather than competing.

Since it is clear that computers can be used to decompress information as well as 
compress it, readers may reasonably ask whether it is sensible to propose that `computing' may be understood as IC. What I believe is a satisfactory resolution of the paradox of ``decompression by compression'' is discussed briefly in Section \ref{paradox} and more fully in \citet{wolff_2000}.

\subsection{PRESENTATION}

The main sections of this article are as follows:

\begin{itemize}

\item Section \ref{preliminaries}, next, is a preparation for two main sections that follow.

\item Section \ref{ML_structures}, which is the first of the two main sections, presents arguments and examples in support of the idea that ML forms and structures can often be 
understood in terms of IC. The possible application of the ICMAUS concepts is 
considered briefly.

\item Section \ref{ML_processes}, the second of the two main sections, describes how a range of ML processes of calculation and inference can be understood in terms of IC, and shows 
with examples and output from the SP61 model how the ICMAUS concepts may 
be applied.

\item Section \ref{discussion_conclusion} briefly considers a range of issues that relate to the ideas that have been described.

\end{itemize}

\section{\bf Preliminaries}\label{preliminaries}

\subsection{REPRESENTATION, SEMANTICS AND PLATONISM}\label{repr_semantics_platonism}

In mathematics and logic, it is often convenient to make a distinction between the representation of an ML entity and the thing itself (the meaning or `semantics' associated with the representation). Thus the concept {\it four} may be represented by `four', `4', `IV', 
`1111' (in unary arithmetic), `100' (in binary arithmetic), and so on.

In mathematical Platonism, mathematical entities ``are not merely formal or quantitative structures imposed by the human mind on natural phenomena, nor are they only mechanically present in phenomena as a brute fact of their concrete being. Rather, they are numinous and transcendent entities, existing independently of both the phenomena they order and the human mind that perceives them.'' \citep[][pp. 95--96]{hersh_1997}. Such ideas are ``invisible, apprehensible by intelligence only, and yet can be discovered to be the formative causes and regulators of all empirical visible objects and processes.'' ({\it ibid.} p. 95).

The view of mathematical (and logical) concepts that has been adopted here is quite different. It is assumed that both the representations of ML concepts and the concepts themselves are `information' in a technical sense (to be described in Section \ref{information_compression}). If it is accepted that brains and computers are valid vehicles for ML concepts and ML processes of calculation and inference, it is difficult to resist the conclusion that ML concepts must exist within those systems in some relatively concrete form---arrays of binary digits in digital computers or, in brains, something like `cell assemblies', patterns of nerve impulses or, perhaps, DNA or RNA sequences. There is no place for ``numinous and transcendent entities'' in a Platonic world outside the brain or computer.

Do G{\"o}del's incompleteness theorems pose a problem for this non-Platonist view? Discussion of that issue would take us well beyond the scope of this article. We can at least say that any such problem is neither more nor less than it would be for any artificial computing system.

\subsection{INFORMATION COMPRESSION}\label{information_compression}

As a preparation for the rest of the article, this subsection presents an outline 
view of IC, picking out what I regard as the key ideas. A much fuller and more technical 
presentation may be found in text books such as \citet{cover_thomas_1991}.

The advent of computers has made familiar the idea that `information' of all kinds 
- speech, writing, paintings, diagrams, music, and so on can be translated into arrays of 
binary digits (`0' and `1'---`bits') stored in computer memory or on disk. 
These arrays may be one-dimensional sequences, or arrays in two or more dimensions.

\subsubsection{\em Randomness and Redundancy}

Arrays or {\it patterns} of information like those just mentioned may be totally {\it random} (something that rarely occurs in practical applications) and we judge the information to be totally lacking in `structure'.\footnote{In this research, the term {\it pattern} has been adopted to mean an {\it array} of one or more {\it symbols}, including a one-dimensional array or {\it sequence} of symbols, an array of symbols in two or more dimensions and also a subset of symbols within a larger array where the symbols in the subset need not be contiguous within the larger array. Although most of the examples in this article are about sequences, the term {\it pattern} will generally be used in this article in preference to the word {\it sequence} for the sake of generality. Notice that an array of symbols can be as small as a single symbol. 

In this research, the term {\it symbol} means an atomic `mark' that can be differentiated in a yes/no manner from other symbols. 

The meanings of {\it pattern}, {\it symbol} and other concepts used in this research are defined in Appendix A of \citet{wolff_1999_comp}.} 

Most of the information we encounter in practical applications is not random but is 
structured in some way. This means that the information contains {\it redundancy}, a technical 
concept with a meaning that is quite close to the everyday meaning of 
`surplus to requirements'. If there is redundancy in a given body of information, $I$,\footnote{Here and throughout this article, the phrase ``a given body of information'' or ``the given body of information'' will be abbreviated as `$I$'.} there is, in effect, repetition of information. And information that is repeated is, in terms of communication, unnecessary.\footnote{Of course, redundancy can be useful in other ways: guarding against loss of information (as in the use of backup copies or mirror disks in computing), in speeding up processing in databases that are distributed over a wide area, or in aiding the correction of errors when information is corrupted (as in error-correcting codes).}

\subsubsection{\em Hartley-Shannon Information Theory and Algorithmic Information Theory}\label{information_theory}

{\sloppy In the Hartley-Shannon concept of information, $I$ contains the maximum possible information for the given alphabet of symbol types if the absolute and conditional probabilities of all the symbol types in $I$ are equal \citep[see][]{cover_thomas_1991}. In this case, $I$ is {\it random} and contains no redundancy. If there are variations amongst the probabilities, $I$ contains less than the maximum possible information for the given symbol types and, for that reason, it is not random and contains redundancy. The redundancy in $I$ is the difference between the information it actually contains and the maximum possible information for the given alphabet of symbol types.}

These concepts of randomness and redundancy are often useful but can be difficult 
to apply in a sensible way if $I$ is small. Alternative but complementary concepts, that do not suffer from this weakness, form part of the foundation of Algorithmic Information Theory (AIT), pioneered by Solomonoff, Kolmogorov, Chaitin and others \citep[see][]{li_vitanyi_1997}.

The key idea in AIT is that $I$ is random if no computer program can be discovered or designed that generates $I$ but is smaller than $I$ in the sense that it contains fewer symbols. If any such program can be discovered or designed, then $I$ is not random and contains redundancy. For all instances of $I$ except very small ones, it is not possible to prove that there is no computer program that will generate $I$ but is smaller than $I$---because failure to find or devise such a program within a given time does not exclude the possibility that one may be found later. Thus, in the great majority of cases, it is not possible to prove that any $I$ is random although the opposite may be proved in particular cases.

\subsubsection{\em Redundancy and the Repetition of Symbols and Patterns}\label{IC_repetition_of_patterns}

Very often, redundancy in $I$ takes the form of repeating patterns, including the 
repetition of single symbols. Consider, for example, a sequence like this: `a a a a a a a a a a b a a a a a a a a a a a a a a a'. In this sequence, the frequency and thus the (absolute) probability of `a' is very much higher than the frequency or probability of `b' so the sequence clearly contains redundancy in terms of the Hartley-Shannon theory. Likewise, the sequence contains redundancy in terms of the AIT concepts of randomness and 
redundancy because it can be compressed into a `program' for generating the sequence 
that looks something like this: `$pr$(a[10], b, a[15])', where `{\it pr}' means `print'. 

In other examples of $I$ containing redundancy, it may be much less easy to see the 
redundancy as repetition of patterns. And in some cases, such as long sequences of the 
decimal values of $\pi$, it may be very hard to see any kind of structure or redundancy at all, although decimal expansions of $\pi$ (apart from very small ones) are known to be highly redundant because they can be generated by a relatively short computer program. Notwithstanding the difficulty of seeing repeating patterns in examples like this, a working hypothesis in this research is that all kinds of redundancy may, ultimately, be understood in terms of the repetition of patterns.

\subsubsection{\em IC and the Unification of Matching Patterns}\label{IC_unification}

IC is achieved by the removal of some at least of the redundancy in $I$, with or 
without the removal of some of the non-redundant information too.

If redundant information is removed from $I$ and all the non-redundant information 
is retained, then in principle and usually in practice it is possible to reconstruct $I$ in its original form with complete fidelity. This is {\it lossless} IC. Clearly, lossless compression is not possible if there is no detectable redundancy in $I$.

If some of the non-redundant information is removed from $I$ together with the removal of redundant information, it is never possible to reconstruct $I$ with complete fidelity. This is {\it lossy} IC. If the degradation of $I$ is not too great, lossy IC may be adopted for some applications because it can be achieved more quickly than lossless IC, or more cheaply, or because it leads to greater compression, or all these things. 

The sobriquet `SP' is sometimes used as an alternative name for ICMAUS 
and associated ideas because IC may be understood as a reduction in redundancy in a body 
of information---thus increasing {\it Simplicity}---whilst preserving as much as possible of the non-redundant descriptive {\it Power}.

The fact that redundancy often manifests itself as repeating patterns (Section \ref{IC_repetition_of_patterns}) suggests a simple method for finding and removing redundancy: look for patterns that repeat in $I$---by systematically comparing or {\it matching} patterns against each other---and then merge or {\it unify} two or more instances of a given pattern to make a single instance.\footnote{Notice that the use of the terms `unify' and `unification' in this article is distinct from the use of those terms in logic. Here, they refer to a simple merging of patterns that match each other, a meaning that is related to but simpler than the meanings of those terms in logic.} All the simpler `standard' techniques for IC (such as PkZip) operate in this way \citep[see, for example,][]{storer_1988}. A corollary of the working hypothesis that all kinds of redundancy may be seen as repeating patterns is that all kinds of technique for IC, including the more elaborate ones, may be seen to work by the unification of matching patterns.

Notice that the concept of {\it counting} the number of instances of a pattern implies the 
matching of patterns against each other to establish that they are instances of the given 
pattern. And the concept of {\it frequency} of a given pattern implies the unification of the 
several instances to form a representative `type' instance of the pattern. Since the concept of probability in this context is a normalised expression of frequencies, it is apparent that there is a strong connection between the concept of probability used in Hartley-Shannon information theory and the notions of matching and unification of patterns. Of course, all of the foregoing also applies to symbols and symbol types.

\subsubsection{\em Three Techniques for IC by the Matching and Unification of Patterns}\label{IC_techniques}

This section describes three techniques for achieving IC by matching and unification of patterns and principles 
that apply in all three cases. The first technique, which we shall call 
{\it chunking-with-codes} is the most basic technique. The other two may be seen to be variants of it.

\begin{itemize}

\item {\it Chunking-with-Codes}. If two or more matching patterns in $I$ are unified to form a single `type' pattern, the result is lossy compression. This is because the location of each of the original patterns is lost. To avoid this loss of non-redundant information, the chunk (pattern) of information may be assigned a relatively short name, `tag', `identifier' or `code'.\footnote{The meaning of the word `code' in this context is `name' or `identifier' and not the same as the meaning of the word in cryptography.} Then copies of this code are placed in each location where copies of the pattern originally appeared, thus preserving the information about the original locations of the given pattern. 

This {\it chunking-with-codes} technique is so widespread and `natural' that we hardly notice it: in this article, the abbreviation `IC' is used as a short `code' for the relatively long pattern `information compression'; a citation like `\citet{storer_1988}' may be seen as a code for the relatively long bibliographic details of the book given in the references section of 
the article; in everyday speaking and writing, names of people, places and so on may all be regarded as relatively short codes for concepts where the full description in each case represents a relatively large body of information. No doubt, readers can think of many other examples.

\item {\it Schema-plus-Correction}. A variant of the basic chunking-with-codes technique is another technique that is often called {\it schema-plus-correction}. In this case, the unified pattern is not a monolithic chunk but is a chunk containing gaps or holes within which a variety of other patterns may appear on different occasions. In this case, the unified pattern is the {\it schema} and the other patterns that fill the gaps are {\it corrections} to or completions of the schema.

In everyday life, a common example is a menu in a restaurant. In this case the 
schema is the basic framework of the menu, e.g., `Starter ... main course ... sweet course ...' and the `corrections' are the dishes chosen to fill the gaps. Another example is any kind of form with fields that will be filled with various `corrections' or completions each time the form is completed. 

\item {\it Run-Length Coding}. If a sequence of symbols contains a pattern that repeats in a sequence of instances that are contiguous, one with the next, then IC can be achieved by reducing the repeating series to one instance with something to mark the repetition or, for fully lossless compression, with something to show the number of repetitions. For example, the pattern `a x y z x y z x y z x y z x y z x y z x y z x y z x y z x y z b' may be reduced to something like `a (x y z)* b' in the case of lossy compression or, for lossless compression, something like `a (x y z)[10] b'.

This kind of {\it run-length coding} can always be expressed using a `recursive' function---one that contains one or more calls to itself, either directly or via calls to other functions (see Section \ref{iteration_and_recursion}).

\end{itemize}

\subsubsection{\em Using Short Codes for Frequent Patterns and Long Codes for Rare Patterns}

Other things being equal, there is a clear advantage in terms of IC if pattern types that occur frequently in $I$ have shorter codes than ones that are rare. A well-known scheme that generates variable-length codes precisely geared to this principle is Huffman coding \citep[see][]{cover_thomas_1991}.

\subsubsection{\em The Complexity of Searching and the Need for Constraints}\label{complexity_constraints}

The idea of matching and unifying patterns is simple enough in itself but can be 
complex to apply, especially if we allow partial matches like the matching of `x z' with the first and third characters in `x y z'. For most patterns, the size of the search space is astronomically large. Given that search spaces are normally so large, it is not normally possible to search them exhaustively. In the vast majority of cases, it is necessary to constrain the search in some way, either eliminating parts of the search space {\it a priori} (e.g., by restricting the search to exact matches) or by applying heuristic techniques or both these things.

In heuristic search, the search is conducted in stages with a narrowing of the search 
space at each stage using some measure of `goodness' to choose where the next stage of 
searching will be concentrated. In the present context, it is appropriate to use some kind of measure of IC as the measure of `goodness' of matches.

\subsection{\sloppy THE PROPOSALS IN RELATION TO ESTABLISHED IDEAS}\label{established_ideas}

\sloppy{The various `isms' in the philosophy and foundations of mathematics and logic---foundationism, logicism, intuitionism (constructivism), formalism, Platonism, neo-Fregeanism, humanism, structuralism and so on---are well described in various sources (see, for example, \citet{potter_2000, hersh_1997, craig_1998, hart_1996, barrow_1992, eves_1990, epstein_carnielli_1989, boolos_jeffrey_1980}) and there is no intention to describe them again here.}

Insofar as the embryonic ideas to be presented constitute any kind of philosophy of 
ML, it seems that they do not fit easily into any of the existing isms or schools of thought about the philosophy and foundations of ML.

\subsubsection{\em Possible Connections}\label{possible_connections}

Devlin's interesting ideas about logic and information \citep{devlin_1991} might be thought to 
be related to the present proposals because a concept of information is fundamental in 
ICMAUS. However, Devlin develops a concept of information that is different from the 
`standard' concepts of Hartley-Shannon information theory and AIT that have been adopted in 
the present work. Perhaps more important is the fact that IC does not have any role in his proposals in contrast to the ICMAUS proposals where IC has a key significance.

Two recent books about mathematics suggest in their titles that mathematics is a 
``science of patterns'' \citep{devlin_1997, resnik_1997}. Given the use of the word 
`pattern' in the present research, one might think that these two books give an account of mathematics that is closely related to the ICMAUS concepts. However, the word 
`pattern' in the present work is merely a generic name for an array of atomic symbols in 
one or more dimensions whereas in \citet{devlin_1997} and \citet{resnik_1997}, the word `pattern' has more abstract meanings.

That said, Devlin's book---which is intended as a `popular' presentation of 
mathematics rather than an academic treatise---gives a hint of a connection with the present work because the term `pattern' is used in a way that seems to be roughly equivalent to `regularity' and there is a connection between regularities and IC. Resnik uses the term `pattern' to mean the same as the abstract concept of `structure' as it has been developed in the {\it structuralist} view of mathematics \citep[see also][]{shapiro_2000}. As with Devlin's book, there is a hint of a connection with the present work because of an apparent link between structuralist concepts and the intuitive concept of structure which itself seems to reflect redundancy.

Despite these possible links with the present work, neither of the books cited 
develops any connection with established concepts of information, redundancy or IC and 
there is no mention of any concept of multiple alignment. Readers may wish to view the present proposals as an extension or development of structuralism. Alternatively, the proposals may be seen as an application of an independent intellectual tradition---{\it information theory}---to concepts in mathematics and logic, with some points of similarity with structuralism.

\subsubsection{\em IC in the Workings of Brains and Nervous Systems}\label{brains_and_ns}

Parsimony has been recognised as an important principle in thinking, science and 
the workings of brains and nervous systems from at least as far back as the fourteenth 
century when William of Occam suggested that, in thinking and theorising, ``entities are not to be multiplied beyond necessity''. Other thinking in this area includes:

\begin{itemize}

\item Ernst Mach and others in the nineteenth century wrote on this theme and, in a similar vein, \citet{zipf_1949} published {\it Human Behaviour and the Principle of Least Effort}.

\item The advent of Hartley-Shannon information theory (ca. 1949) led to a surge of interest in `cognitive economy' and economical encoding in brains and nervous systems \citep[see, for example,][]{attneave_1954, oldfield_1954, garner_1974, von_bekesy_1967}. H. B. Barlow has published extensively on neurophysiology over many years with IC as a recurring theme \citep[see, for example,][]{barlow_1959, barlow_1969, barlow_1997, barlow_2001_bbs}.

\item The idea that pattern recognition may be understood in terms of IC has been 
recognised for some time \citep[see, for example,][]{watanabe_article_1972}.

\item In an extensive programme of research by the author, modelling first-language 
learning by children, MLE principles (see Section \ref{inf_MLE_PR}, below) emerged as a 
unifying theme (see \citet{wolff_1988, wolff_1982} and earlier publications cited there).

\item In everyday life, we understand our world to be populated by discrete, persistent 
`objects', including discrete entities like words. Whatever else they may be, these entities seem to have the role of information `chunks' and, as such, they are a powerful aid to economy in our encoding of the world (Section \ref{IC_techniques}, above). Relevant discussion may be found in \citet{wolff_1993, wolff_1988, wolff_1982} and earlier publications cited there.

\item Many pointers to more recent work relating IC to the workings of brains and 
nervous systems may be found in \citet{chater_1999, chater_1996}.

\end{itemize}

If it is accepted that IC does indeed have a fundamental role in the way brains and 
nervous systems operate, it is reasonable to ask what advantage this may have in natural 
selection.

The most obvious answer is that, given finite resources in brains and nervous 
systems, compression of information enables us to achieve more than would otherwise be 
the case---or use fewer resources for a given level of performance. This would apply to the storage of information and also to the internal transmission of information and the input and output of information.

Although these things are probably important, it seems likely that there is a second 
reason that is at least as important. This is that there is an intimate connection between IC and inductive prediction of the future from the past (described briefly in the next subsection). For any creature, an ability to learn from experience and anticipate positive or negative events in the future must be indispensable for survival.

\subsubsection{\em Inductive Inference, Minimum Length Encoding and Probabilistic Reasoning}\label{inf_MLE_PR}

Largely independent of the research just described, there has been another 
programme of research elucidating the theoretical significance of IC in the abstract problem of inductive inference: probabilistic prediction of future events from patterns of events in the past. This programme of research on `Minimum Length Encoding' (MLE) was pioneered by Solomonoff, Wallace and Boulton, and Rissanen amongst others \citep[see][]{li_vitanyi_1997}.

The connection between IC and inductive inference---which may at first sight seem obscure---lies in the identification of recurrent patterns. As we have seen (Section \ref{IC_repetition_of_patterns}), IC may be achieved by the unification of recurrent patterns. And a repeating pattern like `Spring, Summer, Autumn, Winter' informs out expectations that Summer will follow Spring, Autumn will follow Summer, and so on.

The key idea in MLE is that, in grammar induction and related kinds of processing, one should seek to minimise ($G + E$), where $G$ is the size (in bits) of the `grammar' (or comparable structure) under development and $E$ is the size (in bits) of the raw data when it has been encoded in terms of the grammar. This principle achieves a compromise between grammars that are  very small but inefficient for coding and grammars that can encode data economically but are unduly large \citep{solomonoff_1986}.\footnote{In the SP61 model (that provides all the examples of multiple alignment shown in this article), $G$ has a fixed size. This means that it is only necessary to seek to minimise $E$. In the SP70 model, currently under development, $G$ and $E$ can both vary so both of them must be evaluated in the search for minimum length encoding.}

\subsection{THE ICMAUS FRAMEWORK AND THE SP61 MODEL}\label{icmaus_framework_and_SP61}

The ICMAUS framework is founded on MLE principles, just described. The framework is intended as an abstract model of any natural or artificial system for `computing' or `cognition'. It is envisaged that any such system may be seen to receive `raw' (or `New') data from the environment and to compress it as much as possible by matching it against itself and against already-stored patterns (designated `Old'), unifying those parts that match (thus achieving the effect of `recognition') and storing the parts that do not match (which seems to be a key to `learning'). 

This framework is partially-realised in a software model, SP61, described quite fully in \citet{wolff_2000}. SP61 performs the recognition part of the ICMAUS process but does not attempt learning. Work is currently in progress to extend the model to unsupervised inductive learning. Notwithstanding this limitation of the current model, it is good enough to provide illustrations for many of the examples described below.

Currently, all kinds of information for the SP61 model are expressed as one-dimensional {\it patterns} of atomic symbols (examples will be seen below). It is envisaged that, at some stage, the model will be generalised to patterns in two or more dimensions.

At the core of the ICMAUS framework (and the SP61 model) is a process for finding full matches between patterns or good partial matches. It is a refined version of `dynamic programming' \citep[see, for example,][]{sankoff_kruskall_1983} with advantages compared with standard methods: it can process patterns of arbitrary length, it can find several alternative alignments for any given set of patterns, and the thoroughness of searching can be controlled by parameters.

The matching process is applied recursively so that the system can build up `multiple alignments', with two or more levels, as will be seen below.

Notwithstanding astronomically large search spaces for most multiple alignments, the heuristic techniques used in the SP61 model mean that its computational complexity (time complexity and space complexity) is within acceptable polynomial limits.

To date, application of the ICMAUS framework and the SP61 model include: natural language processing \citet{wolff_2000}; modelling a Post Canonical System (and hence a Universal Turing Machine) \citet{wolff_1999_comp}; and probabilistic and human-like reasoning \citet{wolff_1999_prob, wolff_2001_igpl}.

\subsection{EXAMPLE}\label{parsing_example}

As a typical example of the kind of alignment formed by SP61, Figure \ref{parsing_alignment} shows how the sentence `j o h n r u n s' (supplied to the model as `New') may be aligned with patterns (stored in `Old') that represent grammatical rules. The overall effect of this alignment is an analysis or parsing of the sentence into its constituent parts.

\begin{figure}[!bhpt]
\begin{center}
\begin{BVerbatim}
0       j o h n        r u n s       0
        | | | |        | | | |       
1   N 0 j o h n #N     | | | |       1
    |           |      | | | |       
2 S N           #N V   | | | | #V #S 2
                   |   | | | | |     
3                  V 1 r u n s #V    3
\end{BVerbatim}
\end{center}
\caption{\small Parsing of the sentence `j o h n r u n s' as an alignment amongst patterns representing the sentence and relevant rules in a grammar.}
\label{parsing_alignment}
\normalsize
\end{figure}

In this and other alignments shown below, New is placed at the top of each alignment with patterns from Old below it. Each appearance of a pattern in any alignment is given a line to itself and, below the top line, the order in which patterns appear (from top to bottom of the alignment) is entirely arbitrary.

A good alignment is defined here as one that provides a basis for an economical 
coding of New in terms of the patterns in Old. There is no space here to describe the method of evaluation which is used in SP61. A full description may be found in \citet{wolff_1999_prob}.

\subsection{PARSING IN MATHEMATICS AND LOGIC}\label{maths_logic_parsing}

To bring us a little closer to mathematics and logic, this subsection presents an example showing how expressions in logic and mathematics may be parsed. This, in itself, is not a process of mathematical or logical inference such as addition or multiplication but it is necessary for the interpretation of ML symbols prior to the application of such processes. The example provides a bridge from the `linguistic' example shown above to other examples presented later showing how mathematical and logical structures and processes may be understood in terms of ICMAUS.

Both mathematics and logic make frequent use of `expressions' that need to be analysed into their constituent parts before they can be evaluated. An example from arithmetic is $(((6 + 4)(16 - 4)) - 6)$. In logic one might write a composite proposition like $\sim((p \wedge q) \Rightarrow ((x \wedge y) \vee (s \wedge p)))$, where the letters are simple `atomic propositions' like ``Today is Tuesday'' or an equivalent `atomic formula' like `today(Tuesday)' and the other symbols have their standard meanings in logic (see the key to Figure \ref{logic_grammar}). Any such analysis is essentially a form of parsing in the same sense as was used in Section \ref{parsing_example}. And, like `linguistic' kinds of parsing, it requires an appropriate grammar.

It may be objected that the examples are already parsed because they contain the 
kinds of brackets that may be used to represent the structure of a sequence of symbols. 
The reason that parsing is needed in examples like these is that, until the constituent parts of a sequence are recognised by the system, the sequence is merely a stream of symbols and the role of the brackets in marking structures has not yet been recognised. Brackets are, in any case, not a necessary part of the examples---they might equally well be presented in Polish notation without the use of brackets.

Figure \ref{logic_grammar} shows part of a grammar for expressions in symbolic logic. Readers will see that it contains three main rules (beginning with `F'). The first one creates an `atomic formula' like `today(Tuesday)'. In a fuller grammar there would be rules for the creation of a variety of such atomic formulae. The second rule is a recursive rule covering cases like `proposition implies proposition', where `proposition' can be simple or complex. The third rule is also recursive, covering all forms of `not proposition'. The patterns beginning with `R' show some of the symbols that may be put between two propositions in the second rule in the grammar.

\begin{figure}[!bhpt]
\begin{center}
\begin{tabular}{l}
F $\rightarrow$ a \\
F $\rightarrow$ ( F R F ) \\
F $\rightarrow$ ( $\sim$ F ) \\
R $\rightarrow$ $\Rightarrow$ \\
R $\rightarrow$ $\wedge$ \\
R $\rightarrow$ $\vee$ \\
\\
etc
\end{tabular}
\end{center}
\caption{\small Part of a grammar for symbolic logic. {\it \bf Key:} `a' = `atomic formula', 
`$\sim$' = `not', `$\wedge$' = `and', `$\vee$' = `or', `$\Rightarrow$' = `implies'.}
\label{logic_grammar}
\normalsize
\end{figure}

Figure \ref{logic_grammar_icmaus} shows the same grammar as the one in Figure \ref{logic_grammar} but it is recast into a form that is appropriate for ICMAUS. The main differences are that the grammar in ICMAUS style leaves out the rewrite arrow, it uses numbers to differentiate the three `F' rules and the three `R' rules, and it puts in a `terminating' symbol for each rule and each `call' to a rule. Each `terminating' symbol is formed by copying the initial symbol and adding `\#' to the front of it.

\begin{figure}[!bhpt]
\begin{center}
\begin{tabular}{l}
F 1 a \#F \\
F 2 ( F \#F R \#R F \#F ) \#F \\
F 3 ( $\sim$ F \#F ) \#F \\
R 1 $\Rightarrow$ \#R \\
R 2 $\wedge$ \#R \\
R 3 $\vee$ \#R \\
\end{tabular}
\end{center}
\caption{\small The grammar shown in Figure \ref{logic_grammar} recast in `ICMAUS' style.}
\label{logic_grammar_icmaus}
\normalsize
\end{figure}

\begin{figure}[!bhpt]
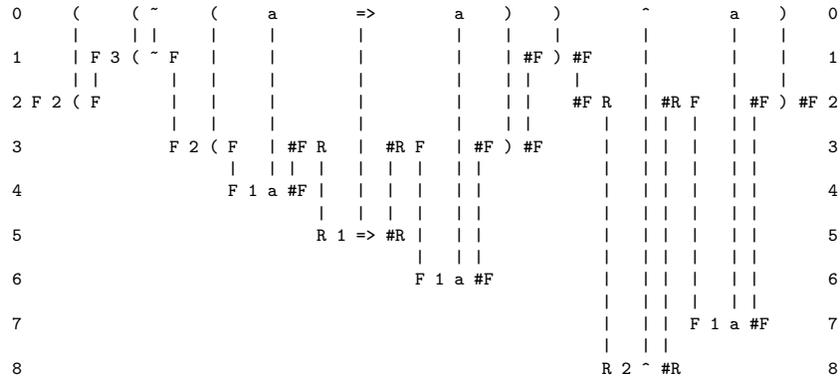

\fontsize{07.00pt}{08.40pt}
\begin{center}
\begin{BVerbatim}
0     (     ( ~     (     a        =>        a    )    )        ^        a    )    0
      |     | |     |     |        |         |    |    |        |        |    |   
1     | F 3 ( ~ F   |     |        |         |    | #F ) #F     |        |    |    1
      | |       |   |     |        |         |    | |    |      |        |    |   
2 F 2 ( F       |   |     |        |         |    | |    #F R   | #R F   | #F ) #F 2
                |   |     |        |         |    | |       |   | |  |   | |      
3               F 2 ( F   | #F R   |  #R F   | #F ) #F      |   | |  |   | |       3
                      |   | |  |   |  |  |   | |            |   | |  |   | |      
4                     F 1 a #F |   |  |  |   | |            |   | |  |   | |       4
                               |   |  |  |   | |            |   | |  |   | |      
5                              R 1 => #R |   | |            |   | |  |   | |       5
                                         |   | |            |   | |  |   | |      
6                                        F 1 a #F           |   | |  |   | |       6
                                                            |   | |  |   | |      
7                                                           |   | |  F 1 a #F      7
                                                            |   | |                
8                                                           R 2 ^ #R               8
\end{BVerbatim}
\end{center}
\caption{\small The best alignment found by SP61 with the expression 
`( ( $\sim$ ( a $\Rightarrow$ a) ) $\wedge$ a )' in New and the patterns shown in Figure \ref{logic_grammar_icmaus} in Old.}
\label{expression_parsing}
\normalsize
\end{figure}

Finally, Figure \ref{expression_parsing} shows the best alignment that was found by SP61 with the pattern `( ( $\sim$ ( a $\Rightarrow$ a) ) $\wedge$ a )' in New and the patterns from Figure \ref{logic_grammar_icmaus} in Old. The alignment may be interpreted as a parsing of the logical expression in terms of the given grammar and, as such, appears to be `correct'.

\subsection{REALISATION OF THE THREE COMPRESSION TECHNIQUES IN ICMAUS}\label{compression_in_ICMAUS}

The two examples of multiple alignment that we have seen so far give us a 
preliminary view of the way in which the three compression techniques described in 
Section \ref{IC_techniques} may be modelled in ICMAUS.

\subsubsection{\em Chunking-with-Codes}

In our first example (Figure \ref{parsing_alignment}), a pattern like `N 0 j o h n \#N' may be seen as a chunk of information (`j o h n') together with its `code' (`N 0 \#N'). In this case, the code has a dual function of providing a unique identifier for the chunk and marking the 
beginning and end of the chunk.

Superficially, the code in this example appears to be only slightly smaller than the 
chunk and there does not seem to be much advantage in terms of IC in recognising the 
chunk. If `j o h n' had been replaced by `j i m', there would appear to be no compression at all. But it must be born in mind that each symbol has its own code and calculations of IC are based on symbol sizes in bits not simple counts of numbers of symbols.

What about patterns like `R 1 $\Rightarrow$ \#R' in our second example where, even allowing 
for varying sizes of symbol codes, it looks as if the code (`R 1 \#R') is larger than what it encodes (`$\Rightarrow$')? This is discussed in Section \ref{SPC}, next.

\subsubsection{\em Schema-Plus-Correction}\label{SPC}

In our first example (Figure \ref{parsing_alignment}), it should be clear that a pattern like `S N \#N V \#V \#S' is taking the role of a schema that may be corrected or completed by the addition of a noun pattern and a verb pattern. Although sentence patterns in a language like English are relatively abstract concepts, they are patterns that repeat very often in the language and thus represent, in compressed form, relatively large amounts of redundancy. In a similar way, patterns like `F 2 ( F \#F R \#R F \#F ) \#F' and `F 3 ( $\sim$ F \#F ) \#F' represent very frequently recurring patterns in the language of logic (`F $\Rightarrow$ ( F R F )' and `F $\Rightarrow$ ( $\sim$ F )' respectively) and thus capture a significant amount of the redundancy in that language.

An answer to the question raised at the end of the last subsection is that, when IC is achieved using schemata like `F 2 ( F \#F R \#R F \#F ) \#F' and `F 3 ( $\sim$ F \#F ) \#F', it is also necessary to create patterns like `R 1 $\Rightarrow$ \#R' that are corrections or completions of the schemata. Such `correction' patterns may not, in themselves, represent any redundancy but they contribute to IC via their role in conjunction with patterns that function as schemata (see also Section \ref{ICMAUS_variables}).

\subsubsection{\em Run-Length Coding}

In ICMAUS, there is no explicit mechanism for representing iterated structures as there is in computer programming (e.g., the {\it while} loop, the {\it for} loop, or the
{\it repeat ... until} loop). The sole mechanism is recursion (in the sense of computer programming) as can be seen in the examples in Sections 
\ref{generation_recognition_numbers} and \ref{maths_logic_parsing}.

The ICMAUS grammar in Figure \ref{logic_grammar_icmaus} exhibits this kind of recursion (in the second and third lines). Another example is described in Section \ref{generation_recognition_numbers}.

\subsection{MATHEMATICS, LOGIC, INFORMATION, REDUNDANCY AND STRUCTURE}

Before we proceed to the two main sections of this article, it may be useful to take a `global' view of ML in terms of concepts that have been introduced in this section.

The concept of `information', as described above, is a highly abstract concept with a very wide scope. Anything and everything of which we may conceive, including mathematics, logic, brains, computers, trees, houses, and so on, may be seen as information. We cannot learn about the world without absorbing information and any inborn knowledge that we may have may also be seen as information.

Information that is totally random---without any detectable redundancy---is totally lacking in structure. All persistent entities, concepts, laws, theorems or regularities, including the concepts of mathematics and logic, represent redundancy or structure in information. Any kind of repeatable calculation or deduction is a manifestation of redundancy in the substance of mathematics and logic. And information compression is the mechanism by which we detect and manipulate these redundancies.

In short, it is difficult to escape from the conclusion that concepts in ML are intimately related at a fundamental level to concepts of information, redundancy, structure and the compression of information.

\subsection{THE ICMAUS FRAMEWORK AND AIT}

AIT (described briefly in Section \ref{information_theory}, above) is closely related to MLE principles. Since the ICMAUS framework is founded on MLE principles, one might ask what distinguishes the ICMAUS proposals from AIT. The key differences are these:

\begin{itemize}

\item Within AIT, the Turing model of `computing' is taken as `given', whereas, within the ICMAUS theory, it may be seen as a special case of the ICMAUS framework \citep{wolff_1999_comp}.

\item The concept of multiple alignment, as it has been developed within the ICMAUS theory, has no counterpart within AIT.

\item The ICMAUS framework provides a perspective on a range of AI issues and on the nature of mathematics and logic which is quite different from anything in AIT.

\end{itemize}

\section{\bf ML Structures and information compression}\label{ML_structures}

The previous sections of this article have prepared the ground for what is really the main 
substance of the article:

\begin{itemize}

\item This section describes the idea that ML notation and ML structures may be seen as 
a set of devices for representing information in a compressed form. 

\item The next section presents the idea that processes of mathematical and logical 
inference may also be understood in terms of IC.

\end{itemize} 

As was noted earlier, the arguments and examples to be presented are not intended 
as any kind of proof of the proposition that ML is based on IC. They are merely suggestive indications of what may turn out to be a fruitful way to view mathematics, logic and related disciplines---an avenue to new insights and understanding.

Within this main section, Section \ref{suggestive_evidence} considers some general features of ML that suggest that there is a close connection between ML and IC. Next, we consider some aspects of number systems that seem to relate to IC. Then in Sections \ref{chunking_with_codes}, \ref{schema_plus_correction} and \ref{run_length_coding}, we consider how the IC principles of chunking-with-codes, schema-plus-correction and run-length coding may be seen to operate in ML forms and structures. Section \ref{modelling_in_ICMAUS} considers some of the ways in which these concepts may be modelled in ICMAUS.

\subsection{SUGGESTIVE EVIDENCE}\label{suggestive_evidence}

Considering, first, the idea that ML notations and structures might be seen as 
devices for representing information in a compressed form, a general indication that 
something like this might be true is provided by three aspects of mathematics discussed in 
the next three subsections.

\subsubsection{\em Mathematics as the ``Language of Science''}

In the context of science as a search for succinct descriptions of the world, it is often 
remarked that mathematics provides a remarkably and `mysteriously' convenient way to 
express scientific truths. Thus \citet{penrose_1989}, who takes a Platonic view of mathematics, 
writes: 

\begin{quotation}
\noindent {\it It is remarkable that} {\bf all} {\it the SUPERB theories of Nature have proved to be 
extraordinarily fertile as sources of mathematical ideas. There is a deep and 
beautiful mystery in this fact: that these superbly accurate theories are also 
extraordinarily fruitful simply as} {\bf mathematics}. 
(pp. 225--226, emphasis as in the original).
\end{quotation} 

\noindent In a similar vein, \citet{barrow_1992} writes: 

\begin{quotation}
\noindent {\it For some mysterious reason mathematics has proved itself a reliable guide 
to the world in which we live and of which we are a part. Mathematics 
works: as a result we have been tempted to equate understanding of the 
world with its mathematical encapsulization. ... Why is the world found to 
be so unerringly mathematical?} (Preface, p. vii).
\end{quotation} 

\noindent This last question is partly answered later: 

\begin{quotation}
\noindent {\it Science is, at root, just the search for compression in the world. ... In short, the world is surprisingly compressible and the success of mathematics in 
describing its workings is a manifestation of that compressibility.} ({\it ibid.}, p. 
247).
\end{quotation}

Thus he recognises that mathematics provides a good medium for expressing scientific 
observations in a compressed form but he stops short of making the connection, which is 
the subject of this article, between mathematical concepts and standard concepts of 
information and IC.

\subsubsection{\em The Nature of Mathematical `Functions'}\label{maths_functions}

A mathematical `function' is a mapping from one set ($X$) to another set ($Y$) and, in 
the `strict' form of a function, each element of $X$ is mapped to at most one element of $Y$ 
\citep[][p. 7]{sudkamp_1988} although it is sometimes useful to consider more `relaxed' kinds of `function' where each element of $X$ may be mapped to more than one element of $Y$. In essence, a function is a table with one or more columns representing the `input' and one or more columns representing the `output'. Thus an element of set $X$ may be one or more cells from one row of the table and the corresponding element of set $Y$ is the complementary set of one or more cells from the same row of the table.

Very simple functions, such as the exclusive-OR function, can be and often are represented literally as a table. This and other examples will be considered later. But, in the vast majority of cases, the definition of a function (the `body' of the function in computer-science jargon) does not look anything like a table. Nevertheless, whatever the definition of a function might look like, all functions behave as if they were `look up' tables that yield an output value or set of values in response to an appropriate input.

The key point here is that, with regard to the more complex kinds of functions that 
do not look like a table, they are almost invariably smaller than the abstract table they 
represent and, in most cases, they are smaller by a large factor.

This feature of the more complex kinds of functions means that they can provide a 
very effective means of compressing information. Consider, for example, a function 
relating the distance travelled by an object to the time since the object was dropped from a high place and allowed to fall freely under the influence of gravity (ignoring the effect of resistance from air). If the function is represented as a table and if there is an entry in the table for each microsecond or nanosecond, then for a fall of 100 metres, say, the table will be very large. By comparison, Newton's formula, $s = (gt^{2}) / 2$ (where $s$ is the distance travelled, $g$ is the acceleration due to gravity and $t$ is the time since the object was allowed to fall), is extremely compact.

This feature of the more complex kinds of function goes a long way to explaining 
why mathematics is such an effective ``language of science''. However, as we shall see, 
there is more to be said than this about IC and the nature of functions.

\subsubsection{\em The Nature of Counting}

Given its close association with the number system, the process of counting things is clearly central in mathematics.

A little reflection shows that it is not possible to count things without some kind of process of recognition. If we are asked to count the number of apples in a bowl of fruit, we need to be able to recognise each apple and distinguish it from other fruits in the bowl---oranges, bananas and so on. Even if we are asked to count the number of `things' in some location, we need to recognise each one as a `thing' and distinguish it from the environment in which it lies.

Recognition is a complex process which can be applied at various levels of abstraction. But notwithstanding this complexity it seems clear that it necessarily involves a process of matching patterns. And, as we noted in Section \ref{IC_unification}, if the things to be counted are all to be assimilated to a single concept (e.g., `apple' or `thing' or other category), this implies the unification of matching patterns.

Given the connection between IC and the matching and unification of patterns (described in Section \ref{IC_unification}), it seems that counting involves a process of compressing information. Given the central importance of counting in mathematics, this is one strand of evidence pointing to the importance of IC in mathematics.

\subsection{CHUNKING-WITH-CODES}\label{chunking_with_codes}

As we saw earlier (Section \ref{IC_techniques}), the idea of giving a name, `identifier', `tag' or `code' to a relatively large chunk of information is a means of compressing information that we use so often and that comes so naturally that we scarcely realise that we are compressing information.

In ML, as in many other areas, names are widespread. We use names for mathematical and logical functions, for sets and members of sets, for numbers, for variables, for matrices and for geometric forms. In computer programming, names are used for functions and subroutines, for methods, classes and objects (in object-oriented design (OOD)), for arrays and variables, for tables and columns (in databases) and, in the bad old days of spaghetti programming, names were used copiously for individual statements or lines where the program could `go to'.

In almost every case, it seems that a name has the role of being a relatively short 
code representing a relatively large chunk of information. As such, most or all of these 
names contribute to the achievement of IC.

In the following subsections we shall consider these ideas in a little more detail.

\subsubsection{\em The Concept of a Discrete Object in ML}

Historically, it seems likely that the concept of number originated with the need to 
keep a tally on concrete objects like sheep, goats, oranges etc. Arguably, the earliest 
number system was the system of unary natural numbers where $1 = 01, 2 = 011, 3 = 
0111$ etc and other number systems were invented or discovered later as a generalization of that simple kind of counting.\footnote{A concept of zero was, historically, a relatively late development.}

Whether or not that historical speculation is accepted, most people will accept that 
there is a connection between the concept of natural numbers, especially unary numbers, 
and the concept of discrete objects in the world. The latter do not depend on the former (because animals---and people without any knowledge of numbers---can recognise discrete objects). But each stroke in a unary number is itself recognised as a discrete entity and, even if we are using numbers with bases higher than one, we know that it is often convenient to make a one-to-one association between strokes in the underlying unary number and objects in the real world.

In a similar way, the concept of a set, with elements of the set being discrete entities 
(some of which may themselves be sets), provides a link with the way we see the world as 
composed of discrete objects.

Thus the concept of a discrete object has a significant role in both the concept of 
number and of a set. And the concepts of number and set are both key parts of ML. Bearing in mind that seeing the world in terms of discrete objects is almost certainly a reflection of the way our perceptions and concepts are governed by IC principles (Section 
\ref{brains_and_ns}), it is apparent that IC plays a significant part in the foundations of ML. 

\subsubsection{\em Chunking-with-Codes in Number Systems?}

With regard to names or codes as relatively short identifiers for chunks, let us begin 
with a type of name which, at first sight, does not appear to have a role in IC. It can be argued that the names (words or numerals) that we give to numbers merely serve to 
differentiate one number from another and are {\it not} codes for relatively large chunks of 
information. We may argue that this is true even for names like `one billion' or `$10^9$' that represent relatively large numbers and would take a lot of space if they were written out using unary arithmetic or even binary arithmetic. The reason we may give is that the base of a number (unary, binary etc) is irrelevant to the calculation of the amount of information it represents so that any number in any base is totally equivalent to the same number in another base in terms of information, apart from errors arising from any rounding of numbers that we may choose to make.

This line of reasoning is valid for conversions between two numbers both of which 
have a base of 2 or more, but it appears to be wrong for conversions between numbers 
where one of the two numbers is unary. 

Although unary numbers are conventionally considered to have a base of 1, there is 
a sense in which they have no base at all. Each stroke in a unary number seems to be a 
primitive representation of a discrete entity and the stroke cannot be analysed into anything 
simpler.

By contrast, any number that uses a base greater than $1$ introduces a new principle: 
the idea that groups of symbols may be given a relatively short name or code, and that this may be done recursively through any number of levels. In the decimal system, for example, $8$ is a relatively short name for the unary number `$0 1 1 1 1 1 1 1 1$'. At this point, a sceptic might object that, to accommodate decimal digits from $0$ to $9$, each one must be assigned at least $\lceil log_2 10\rceil = 4$ bits of space in a computer memory and that this is equivalent to the amount of information in the corresponding unary number.

In the case of the decimal digits from $0$ to $9$, this objection may carry some force. 
But it appears not to be valid when we consider numbers of $10$ or more. The decimal 
number $100$, for example, uses digits within the range $0$ to $9$ and so, on the assumption that we are not needing to represent anything other than numbers, we could accommodate the three digits in `$100$' with only $3 \times 4 = 12$ bits (or, if we started counting from $0$ and thus finished counting at $99$, we could accommodate the number with only $2 \times 4 = 8$ bits). This is very much less than the number of strokes in the corresponding unary number and thus, {\it prima facie}, seems to be a clear case of IC.

In Section \ref{IC_unification} we saw that lossless compression of $I$ is only possible if we can detect redundancy in $I$. Where is the redundancy in unary arithmetic that might allow us to use the chunking-with-codes technique as suggested above? The answer seems to lie in the apparent need to provide a symbol (normally `$0$') to represent zero in the unary number system.\footnote{Apart from the need to represent the concept of zero, a system containing nothing but unary numbers would need some means of marking the beginnings and ends of numbers. This can be achieved if every number starts with `0'.} If the symbols `$0$' and `$1$' are both provided in the unary number system then it is clear that in all the many unary numbers where the number of `$1$'s is large, the probability of `$1$' is very much higher than `$0$' and so, in accordance with Hartley-Shannon information theory, these numbers generally contain considerable amounts of redundancy. Thus the conversion of any unary number (except very small ones) into the equivalent number using a base of $2$ or more seems to be a genuine example of chunking-with-codes and seems to achieve genuine IC.

\subsubsection{\em Grammars and Re-Write Rules}

CF-PSGs comprising rewrite rules like `S $\rightarrow$ NP VP' are used in theoretical linguistics although, since Chomsky's monograph on {\it Syntactic Structures} \citep{chomsky_1957}, it is well known that such grammars are not adequate to handle the full complexity of natural languages. Essentially the same kind of grammar and rules are known as Bachus-Naur form (BNF) and are used widely in computer science, especially in specifying the syntax of programming languages. Again, the same kind of grammar and rewrite rules are used in studies of logic \citep[see, for example,][]{plaisted_1993}.

In many examples of rewrite rules, but not all, there are fewer symbols to the left of 
the rewrite arrow than there are on the right. Very often, there is only one symbol before the rewrite arrow and several to the right.

In cases like these, it seems reasonable to interpret the symbol or symbols to the left 
of the rewrite arrow as a relatively short code for the symbols to the right of the arrow. And in most cases, it is clear that these codes do indeed have a role in IC: a grammar is, in effect, a compressed representation of the natural or artificial language that it generates.

In general, a grammar represents lossy compression of any one sample of the 
language. This is because, in most cases, a grammar does not encode specific sentences or 
computer programs. Nevertheless, the rules of a grammar normally correspond to repeating 
patterns found in samples of the language and thus represent a distillation of some at least of the redundancy to be found in the language. The patterns that repeat in a sample of language can be relatively concrete patterns such as words, or they can be more abstract patterns at `higher' levels, like phrases, clauses and sentences.

\subsubsection{\em Named Functions and Operators}\label{named_functions}

In the course of writing a computer program, it often happens that the programmer 
encounters sections of the program in two or more different places within the program 
that are the same as each other or very similar. In cases like this, it is normal practice to reduce the several sections to a single function, to give the function a name and then to use that name as a `call' to the function in each of the places where one of the original sections was found. This kind of procedure provides us with an illustration of the chunking-with-codes technique and its use in the service of IC. It is so convenient to be able to reuse functions like `date' or `time' under Unix by just giving their names, and it would be so tedious if we had to write out in full the program statements for these functions every time we wanted to use them.

Of course, a named function can itself contain one or more named functions and this 
gives rise to the familiar hierarchies of function calls in computer programs and other areas 
of ML.

Notions of `method' (in OOD) or `subroutine' provide similar examples. An `operator' like `$+$' or `$\times$' is essentially the same except that it is normally supplied by the programming system (although systems like Prolog allow users to define their own operators).

Very often, two or more sections of a program are the same except that different 
numbers or data structures are processed. In this case, a named function can be formed but the variable part is supplied as an `argument' to the function, e.g., `$4$' in `$square\_root(4)$'. More will be said below about this kind of function.

In practice, things are not quite as simple as just described. In conventional 
computers, there is a run-time `overhead' arising from the creation of a named function and, if run-time considerations are pressing, a programmer may deliberately choose to leave repeating sections as they are. In other situations, it can be convenient for reasons of readability or layout of the program to create a named function even though there is only one instance of the relevant section of the program. Notwithstanding examples like these, it remains true that a major reason for the use of named functions is to avoid undue redundancy in computer programs.

Very much the same things can be said about the concept of `function' in 
mathematics. Repeated chunks may not always be converted into a named function and 
there may be occasions where a named function is recognised even though it is used only 
once. But in the great majority of cases, the use of a named function can be seen to be a way of avoiding undue redundancy in mathematical text.

Systems like Prolog that are designed for `programming in logic' provide similar structures. A Horn clause in Prolog with both a `head' and a `body' is very much like a function in procedural programming languages. The head serves as an identifier or code for the clause and the body is the chunk of information which is identified. Normally, the body is significantly larger than the head and, normally, the body is a pattern that would be repeated in two or more parts of the program if it had not been recast as a named chunk in the form of a Horn clause. As with functions, there may be `arguments' or `parameters' to allow for variations in the data being processed.

\subsection{SCHEMA-PLUS-CORRECTION}\label{schema_plus_correction}

This subsection considers some constructs in ML that seem to be examples of the 
IC principle of schema-plus-correction.

\subsubsection{\em Structures Containing Variables}\label{ICMAUS_variables}

The concept of a `variable' is widespread in mathematics, logic, computer 
programming and theoretical linguistics. In all its uses it is understood as some kind of 
`receptacle' that may receive a `value' by `assignment' or by `instantiation'. In many 
systems, each variable has a `type' with a `type definition' that defines the range of values that the variable may receive.

In the majority of cases, a variable has a name, although unnamed variables are used 
sometimes in systems like Prolog. If the size of the value of a named variable is relatively large compared with the name (and values in Prolog can be very large, complex structures), then, in some cases, the name may be understood as a `code' and the value may be seen as a corresponding chunk so that the two together may be seen to make a contribution to IC via the chunking-with-codes technique. But in many cases, the value and the name are similar in size or the value is smaller than the name. Surely, in cases like these, or in cases where the variable has no name, there cannot be any contribution to IC?

The answer to this question seems to be much as was described in Section \ref{SPC}, above. A structure that contains variables may usually be seen as a schema, the variables may be seen as `holes' or `slots' in the schema and the values may be seen as `corrections' or `completions' of the schema. In cases like these, the schema may often be seen as a unification of two or more patterns with a corresponding gain in IC. In short, the contribution that variables can make to IC is not normally a function of the variables in themselves but because they form part of a larger structure---a schema---that contributes to IC. Some examples are described in the subsections that follow.

\subsubsection{\em Functions with Arguments or Parameters}

As we noted in Section \ref{named_functions}, a function need not be a monolithic chunk of information: most functions allow for variations in the data to which they will be applied by providing `arguments' or `parameters' that may take different values on different invocations of the function.

This kind of function may be seen as a schema with holes or slots (the arguments or parameters) that may be filled with a variety of `corrections' or `completions' on different occasions. Provided the function is called two or more times in an ML text, it can be seen as a means of avoiding undue redundancy in that text.

\subsubsection{\em Classes, Objects and Inheritance of Attributes}\label{classes_objects_inheritance}

In our everyday talking, writing and thinking, we make frequent use of general 
classes like `person', `house', `furniture' and so on. Classes may include other classes, 
recursively through any number of levels. The class `person' for example, may include 
subclasses like `man' or `woman' and it may be included in higher level classes like `living thing'.

Although it is difficult to be sure how our brains actually use these kinds of 
concepts, it seems likely that they provide a powerful means for compressing information. 
Certainly, this would be very much in keeping with the variety of other evidence of 
economy in the functioning of brains and nervous systems that we noted earlier (Section \ref{brains_and_ns}).

The trick is to store all the attributes of living things in general in our memory of 
the class `living thing', to store all the {\it additional} things we know about people in general in our memory for `person', and likewise for the classes `man' and `woman'. Then if we recognise something as a man, we can infer that he has all the attributes of `man' and of `person' and of `living thing'. Without this device, information about the attributes of `man', `person' and `living thing' would have to be repeated in the mental record for each of the men that we may know. If, like most people, we know at least a little about each of a large number of different men, the redundancy would be huge.

Starting with Simula (invented by Dahl, Myhrhaug and Nygaard \citep[see, for example,][]{birtwistle_1973}), these ideas have been incorporated in a variety of 
`object-oriented' programming languages. Using such a language, it is possible to define classes like those above and, usually for the lowest level classes, we can create any number of specific objects, each one with its own specific attributes (e.g., `name', `address', 
`telephone number' etc for an object of the class `man'). Any one object {\it inherits} all the `attributes' of its lowest level class and, following the hierarchy upwards, it inherits the attributes of all the higher-level classes of the object. In this programming context, `attributes' are data structures and `methods' (similar to functions). Each object contains its own data structures following the patterns defined in the relevant classes but it can inherit its methods by `borrowing' from the class description as and when they are needed.

This is not the place to rehearse the many benefits of this style of programming. The 
key point here is that it makes it easy to design programs with a minimum of redundancy. 
In terms of the concept of schema-plus-correction:

\begin{itemize}

\item A class at any level is, in effect a schema describing attributes which are found in every member of the class.

\item Any subclass of a given class provides full or partial `corrections' or `completions' of its parent class.

\item The details of any specific object (e.g., the name, address etc of a man) are, 
likewise, corrections or completions of the lowest level class in which the object 
belongs.

\end{itemize}

Notwithstanding the `power' of the concepts of classes, objects and inheritance of 
attributes they seem not to have been taken up in any significant way in mathematics, logic or theoretical linguistics. That said, there have been some moves to incorporate OO principles in versions of Prolog and, as we shall see below, some similarities may be traced between these ideas and others that are prominent elsewhere in ML.

\subsubsection{\em Sets}

Perhaps the closest relative of the OO concepts just described in mathematics and 
logic is the concept of a set.

A set can be defined extensionally---by listing the elements in the set (e.g., $\{sun, 
wind, rain\}$)---or it may be defined intensionally---by listing the defining characteristics of the class (e.g., \{$n | n$ is prime\}). An intensionally-defined set is similar to a class as described in Section \ref{classes_objects_inheritance} and, when one set includes one or more other sets recursively through an arbitrary number of levels, the structure is similar to a class hierarchy as described above. But, to my knowledge, little or no attempt has been made to develop concepts like inheritance of attributes or the instantiation of a class by one or more objects.

With regard to IC, an intensionally defined set that describes two or more 
elements may be regarded as a compressed representation of the relevant attributes of those elements.

\subsubsection{\em Propositions in Logic with the Universal Quantifier}

Here are two examples of propositions expressed in English and in logical notation:

\begin{itemize}

\item \raggedright ``All humans are mortal.''

\indent $\forall x: human(x) \Rightarrow mortal(x)$.

\item ``For any integer which is a member of the set of natural numbers ($N$) and which is 
even, it is also true that the given integer with 1 added is odd.''

\indent $\forall i\ \epsilon\ N \cdot is\_even(i) \Rightarrow is\_odd(i + 1)$.

\end{itemize}

Each example is a proposition containing a universal quantifier and a variable (all three 
instances of $x$ refer to a single variable and likewise for $i$). The first statement is a 
generalization about what it is to be human and the second one is a generalization about 
natural numbers. Both of them represent repeating (redundant) patterns. Each may be 
regarded as a schema rather than a chunk because, notwithstanding its status as `bound', the variable in each example can receive a variety of values.

\subsection{RUN-LENGTH CODING}\label{run_length_coding}

When something repeats in a series of instances, each one contiguous with the next, we find it very natural to find some way to reduce the repetition to a single instance with something to mark the repetition. As noted previously, lossless compression is achieved if the number of 
repetitions is marked; otherwise the IC is lossy.

This section identifies some examples of run-length coding and further examples will be described in Section \ref{ML_processes}.

\subsubsection{\em Iteration and Recursion in Computer Programming}\label{iteration_and_recursion}

In computer programming, repetition of a part of a program can be shown using 
forms like a {\it while} loop or a {\it for} loop or a {\it repeat ... until} loop. Generally speaking, there is some change in the value of one or more variables on each iteration of the block of program. So this kind of repetition may be regarded as a combination of run-length coding and schema-plus-correction.

As we have seen earlier (Section \ref{IC_techniques}), repetition can also be expressed with recursion (in the sense understood in computer programming). Here is an example (in the C programming language) of the factorial function:

\begin{center}
\begin{BVerbatim}
int factorial(int x)
{
     if (x == 1) return(1) ;
     return(x * factorial(x - 1)) ;
}
\end{BVerbatim}
\end{center}

\subsubsection{\em Repetition in Mathematics}

Mathematics uses a variety of ways to show repetition in a succinct manner:

\begin{itemize}

\item Multiplication is a shorthand for repeated addition.

\item The power notation (e.g., $10^9$) is a shorthand for repeated multiplication.

\item A factorial (e.g., $10!$) is a shorthand for repeated multiplication and subtraction.

\item The bounded summation notation (`$\sum$') and the bounded power notation (`$\prod$') are 
shorthands for repeated addition and repeated multiplication, respectively. In both 
cases, there is normally a change in the value of a variable on each iteration, so these 
devices may be seen as a combination of run-length coding and schema-plus-correction.

\item As in our previous examples, repetition of a function may be achieved by a direct 
or indirect `recursive' call to the function.

\end{itemize}

As we noted previously (Section \ref{IC_techniques}), the use of shorthands is so automatic and `natural' that it is easy to overlook the fact that they are examples of information compression.

\subsection{MODELLING THESE CONCEPTS IN ICMAUS}\label{modelling_in_ICMAUS}

If it is accepted that our three basic techniques for IC can be modelled in ICMAUS 
(as outlined in Section \ref{compression_in_ICMAUS}) and if it is accepted that many forms and structures in ML may be understood in terms of IC (as described in the preceding parts of this main section), then one might expect that many of the forms and structures of ML could be modelled in ICMAUS. The majority of examples will be seen in Section \ref{ML_processes}. As a preparation for these later examples, this subsection focuses on some basic structures in ML and how they may be modelled in ICMAUS.

\subsubsection{\em Variable, Value and Type Definition}\label{variable_value_type_definition}

The basic idea of a named variable can be represented quite simply in ICMAUS by 
a pair of symbols that are contiguous within a pattern. In principle, any pair of symbols 
may be seen to function as a variable but in most examples we shall be considering the left 
and right symbols will correspond to the initial and terminal symbols of patterns. Thus, in 
Figure \ref{parsing_alignment}, each of the pairs of symbols `N \#N' and `V \#V' in the pattern `S N \#N V \#V \#S' may be seen to function as variables with the first symbol of each pair acting as a name for the variable and the second serving to mark the end of whatever `value' goes between the two symbols. These `values' are patterns representing a noun in the case of `N \#N' and a 
verb between the symbols `V' and `\#V'. In a similar way in Figure \ref{expression_parsing}, each of the pairs of symbols `F \#F' and `R \#R' in the pattern `F 2 ( F \#F R \#R F \#F ) \#F' function as a variable that can take other patterns in the grammar as values.

In these examples, each variable may be seen to have a `type definition' in the form 
of the set of patterns that has initial and terminal symbols that match the left and right 
symbols of the variable. For the variable `N \#N', the type definition is the two `noun' 
patterns that can fit in that slot and for the variable `V \#V' the type definition is the two 
patterns for verbs. Of course, in a realistic system, those two types would have many more 
members.

Although the examples shown seem to capture the essentials of the concepts of 
`variable', `value' and `type definition', they lack two of the features that are commonly 
found in computer systems:

\begin{itemize}

\item {\it Scoping}: In most computer systems, a variable can have a `scope' meaning a part 
of a program within which repeated instances of a variable name refer to a single 
variable. 

\item {\it Assignment of a type to a variable}: In many computer systems, any variable can be 
assigned a type that is independent of its name. This contrasts with the examples 
given above where the name of the variable is also the name of its type definition.

\end{itemize}

It seems likely that these two ideas can be modelled in ICMAUS but not in the SP61 
model in its current state of development. Discussion of how these things might be done in 
ICMAUS in the future would take us too far afield for this article.

\subsubsection{\em Objects, Classes and Inheritance of Attributes}

Although notions of {\it objects}, {\it classes} and {\it inheritance of attributes} seem not to have been adopted in `mainstream' mathematics or logic, they have become prominent in 
OOD. The example presented here is intended to show the essentials of how these concepts can be accommodated in ICMAUS. As elsewhere in this article, a simple example has been chosen 
for the sake of clarity and to save space but readers should not infer that it represents the 
limits of what the SP61 model can do. The scaling properties of the model are good \citep[see][]{wolff_1999_prob} and more complicated examples are well within its reach.

Figure \ref{animal_taxonomy} shows a set of patterns that represent, in highly simplified form, part of the taxonomic hierarchy of animals recognised by zoologists. The seventh pattern, for example, represents the class `bat' and describes its attributes briefly as `flies echo\_sounder'. Of course, in a realistic example, this description of the attributes of a bat would be much fuller. 

After the symbol `bat' in the pattern for the class `bat' is the name of the class (`mammal') that is immediately above `bat' in the taxonomy of animals. References like this from one class to another means that the patterns capture the (simplified) taxonomic hierarchy, from `vertebrate' and `invertebrate' at the most abstract level down to `Fido' and `Tibbs' at the most concrete level.

\begin{figure}[!bhpt]
\begin{center}
\begin{tabular}{l}
vertebrate backbone \#v \\
invertebrate no\_backbone \#iv \\
mammal vertebrate \#v furry warm\_blood milk \#m \\
bird vertebrate \#v feathers wings \#bd \\
cat mammal \#m retractile\_claws \#c \\
dog mammal \#m barks \#d \\
bat mammal \#m flies echo\_sounder \#bt \\
inst Fido dog \#d \#i \\
inst Tibbs cat \#ct \#i \\
\end{tabular}
\end{center}
\caption{\small A set of patterns representing part of the taxonomy of animals in highly simplified form.  Frequencies of occurrence are assigned to the patterns in accordance with what one would naturally find: `animals' are more frequent than `vertebrates' which are more frequent than `mammals', which are more frequent than `dogs'.}
\label{animal_taxonomy}
\normalsize
\end{figure}

\begin{figure}[!bhpt]
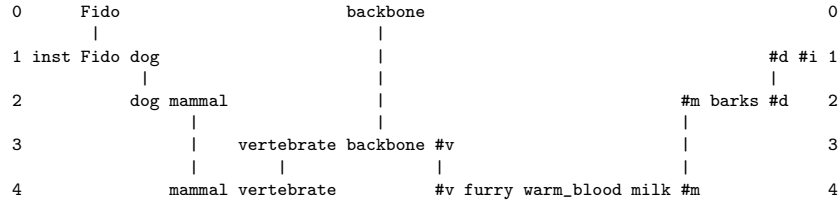

\fontsize{07.00pt}{08.40pt}
\begin{center}
\begin{BVerbatim}
0      Fido                       backbone                                         0
        |                            |                                            
1 inst Fido dog                      |                                       #d #i 1
             |                       |                                       |    
2           dog mammal               |                              #m barks #d    2
                  |                  |                              |             
3                 |    vertebrate backbone #v                       |              3
                  |        |               |                        |             
4               mammal vertebrate          #v furry warm_blood milk #m             4
\end{BVerbatim}
\end{center}
\caption{\small The best alignment found by SP61 with `Fido backbone' in New and the patterns from Figure \ref{animal_taxonomy} in Old.}
\label{fido_alignment}
\normalsize
\end{figure}

Figure \ref{fido_alignment} shows the best alignment found by SP61 with the pattern `Fido backbone' in New (describing something called `Fido' that has a backbone) and the patterns from Figure \ref{animal_taxonomy} in Old. In effect, the alignment has shown that Fido belongs in the class `dog', that this is a `mammal' and that the class `mammal' belongs in the class `vertebrate'. All the attributes of each class are brought into the alignment and, together, they 
provide a relatively full description of the attributes of Fido. This is the essence of the 
concept of `inheritance of attributes'.

\section{\bf ML processes and information compression}\label{ML_processes}

The previous section presented a range of examples showing that many notations 
and forms in ML may be seen as devices for representation of ML concepts in a succinct 
form. This section presents a range of examples in support of the idea that, very often, the 
computational or inferential {\it processes} in ML may also be understood in terms of IC. More 
specifically, it is suggested that, very often, they may be understood in terms of ICMAUS.

As previously noted, an argument has already been presented \citep{wolff_1999_comp} that 
the concept of `computing', as defined by the structure and functioning of the UTM and 
PCS models, may be understood in terms of ICMAUS. The arguments and examples here 
may be seen as a development of the theme with other examples in mathematics, logic and 
related disciplines.

\subsection{SETS}

This subsection considers briefly how some of the basic features of sets may be 
understood in terms of the alignment and unification of patterns.

\subsubsection{\em Creating a Set from a Bag or Multiset}

The concept of `New' as it has been described in this article allows two or more 
patterns within New to be identical and, in that respect, New is like a `bag' or `multiset' in 
logic.

\sloppy If New were to contain a range of one-symbol patterns like this: \{(A)(B)(C)(D)(A)(D)(B)(A)(C)(C)(A)(C)\} (or if New were to contain multi-symbol patterns that repeated exactly in a similar way) it is not hard to see how, with Old initially empty, the alignment and unification of patterns in a learning process like that outlined in Section \ref{icmaus_framework_and_SP61} would reduce New to a set of patterns in Old like this: \{({\bf A})({\bf B})({\bf C})({\bf D})\}.

Thus, the process of converting a bag into the corresponding set containing a single 
example of each type of element in the bag may be seen as a process of information 
compression by the alignment and unification of patterns.

\subsubsection{\em Union and Intersection of Sets}

If the set \{(B)(C)(D)(F)(G)\} were in New and the set \{(A)(B)(C)(E)(F)\} were 
the entire contents of Old, it is not hard to see that, in much the same way as just described, the result would be: \{(A)({\bf B})({\bf C})(D)(E)({\bf F})(G)\}. This set is the {\it union} of the first two sets and their intersection is the one-symbol patterns show in bold type---the set \{({\bf B})({\bf C})({\bf F})\}.

In general, we can see that the union and intersection of two sets may be seen as IC 
by alignment and unification amongst the elements of the sets.

\subsubsection{\em Intensionally-Defined Sets}

An intensionally-defined set such as \{$x | x > 1$\} may be regarded as a 
representation of a set in compressed form. Unless the set has only a few members (or is 
empty), this representation is clearly more economical than the extensional definition 
that itemises every element of the set. The process of deriving an intensional definition 
from the corresponding extensional definition may be seen as the alignment and unification 
of the defining attributes in the elements of the set.

\subsection{GENERATION AND RECOGNITION OF NUMBERS}\label{generation_recognition_numbers}

The simplest system of numbers, in terms of numbers of symbol types, is the unary system mentioned earlier. This kind of system can be defined using a PCS \citep{post_1943}, like this:

\begin{itemize}

\item {\em Alphabet}: 0, 1.

\item {\em Axiom}: 0.

\item {\em Production}: If any string `\$' is a number, then so is the string `\$ 1'. This can be expressed with the production: $\$ \rightarrow \$\ 1$
\end{itemize}

Here, the production is a rewrite rule of the kind discussed earlier and `\$' represents a 
variable. With this `successor function', it should be clear that, if the variable is initialised with the value of the axiom ($0$), then the system will generate, recursively, the infinite set of unary natural numbers, $01, 011, 0111$ etc. Slightly less clear is the fact that the system may be run `backwards' so that it can recognise any unary number as a unary number and thus distinguish it from all other strings which are not unary numbers.

The axiom and production can be modelled with ICMAUS-style patterns shown in 
Figure \ref{unary_icmaus_patterns}. In these patterns, the symbols `a' and `b' are needed for the scoring system used in SP61 and may be ignored. Notice that symbol `\$' in Figure \ref{unary_icmaus_patterns} is not a variable. It is an atomic symbol with exactly the same status as the other symbols in the figure. In functional terms (as described in Section \ref{variable_value_type_definition}, above), the `variable' in Figure \ref{unary_icmaus_patterns} is the pair of symbols `\$ \#\$'.

The recursive nature of the production shown above is represented in Figure \ref{unary_icmaus_patterns} in 
the last pattern shown in Old by the fact that the initial and terminal symbols of the pattern 
match the two symbols `\$ \#\$' which, together, constitute the `variable' in the body of the 
pattern. This pattern may be read as ``A unary number is any unary number (in the `variable' 
in the body of the pattern) followed by `1'''. There is no need for the rewrite arrow which 
appears in the production, above.

\begin{figure}[!bhpt]
\begin{center}
\begin{tabular}{l}
{\it New}\\
\\
$0$\\
\\
{\it Old}\\
\\
\$ a 0 \#\$\\
\$ b \$ \#\$ 1 \#\$\\
\end{tabular}
\end{center}
\caption{\small Patterns for processing by SP61 to model a PCS to generate unary numbers.}
\label{unary_icmaus_patterns}
\normalsize
\end{figure}

\subsubsection{\em Generation of Unary Numbers}

Figure \ref{unary_production_alignment} shows one of many possible `good' alignments between the one-symbol 
pattern in New in Figure \ref{unary_icmaus_patterns} and the patterns in Old. The symbol `0' in New aligns with the matching symbol in the `axiom' in Old. And then that pattern matches the `production' pattern which then matches its initial and terminal symbols to its `variable' symbols recursively.

\begin{figure}[!bhpt]
\begin{center}
\begin{BVerbatim}
0                     0                        0
                      |                       
1                 $ a 0 #$                     1
                  |     |                     
2             $ b $     #$ 1 #$                2
              |              |                
3         $ b $              #$ 1 #$           3
          |                       |           
4     $ b $                       #$ 1 #$      4
      |                                |      
5 $ b $                                #$ 1 #$ 5
\end{BVerbatim}
\end{center}
\caption{\small One of the many good alignments found by SP61 with the one-symbol pattern `0' 
in New as shown in Figure \ref{unary_icmaus_patterns} and with patterns in Old as shown in that figure.}
\label{unary_production_alignment}
\normalsize
\end{figure}

All the symbols in the alignment other than the `0' and `1' symbols may be regarded 
as `service' symbols and may be ignored. If the alignment is unified to form a single 
pattern, it contains the sequence `$0 1 1 1 1$' which is the unary number for 4. In short, the 
alignment may be seen as the generation of that number.\footnote{Readers may, very reasonably, ask ``How is it possible for a relatively large pattern like `$0 1 1 1 1$' to be the outcome of the compression of a relatively small pattern like `$0$'?'' Discussion of the paradoxical way in which ICMAUS allows IC to have the effect of creating something bigger than the original pattern is discussed briefly in Section \ref{paradox}, below, and more fully in \citet{wolff_2000}.}

\subsubsection{\em Recognition of a Unary Number}

Figure \ref{unary_recognition_alignment} shows the best alignment that is found by SP61 with the pattern `$0 1 1 
1 1$' in New and patterns in Old as shown in Figure \ref{unary_icmaus_patterns}. This alignment may be interpreted as the recognition of the binary number for 4 in terms of patterns in Old as shown in Figure \ref{unary_icmaus_patterns}.

\begin{figure}[!bhpt]
\begin{center}
\begin{BVerbatim}
0                     0    1    1    1    1    0
                      |    |    |    |    |   
1                 $ a 0 #$ |    |    |    |    1
                  |     |  |    |    |    |   
2             $ b $     #$ 1 #$ |    |    |    2
              |              |  |    |    |   
3         $ b $              #$ 1 #$ |    |    3
          |                       |  |    |   
4     $ b $                       #$ 1 #$ |    4
      |                                |  |   
5 $ b $                                #$ 1 #$ 5
\end{BVerbatim}
\end{center}
\caption{\small The best alignment found by SP61 with `0 1 1 1 1' in New and patterns in Old as shown in Figure \ref{unary_icmaus_patterns}.}
\label{unary_recognition_alignment}
\normalsize
\end{figure}

\subsection{EXECUTION OF SIMPLE FUNCTIONS}

The concept of a `function', that we have already considered as an example of 
chunking-with-codes (Sections \ref{IC_techniques} and \ref{named_functions}) and as a mediator of IC (Section \ref{maths_functions}), is one of the most distinctive features of ML, especially in mathematics and computer programming.

It is often convenient to think of a function as some kind of `black box' that can 
receive some kind of `input' and, after some processing, will deliver some `output'. As we 
noted in Section \ref{maths_functions}, a function can, in principle, be regarded as a lookup table and, in simple cases, functions can often be implemented in that way. But in more complex cases, a simple lookup table would not be appropriate, either because the size of the table would be finite but too large to be practical or because the table would be infinitely large.

The remainder of this subsection describes with an example how the execution of 
simple functions may be modelled in ICMAUS and the next subsection considers how the 
execution of the more complicated kinds of functions may be understood within ICMAUS.

\subsubsection{\em A Simple Function as a Look-up Table: A One-Bit Adder}

Where a simple function is structured as a table, it can be processed by searching 
for a row in the table where the input matches the input side of the row and then reading off 
the corresponding output.

In ICMAUS, a table may be represented with a set of patterns like those shown in 
Figure \ref{half_adder_definition}. These patterns define the function for the addition of two one-bit numbers in binary arithmetic. The two digits immediately following `A' in each row are the two one-bit numbers to be added together, the bit following `S' is the sum of the two input bits, and the bit following `C' is the {\it carry-out} bit, to be carried to the next column. Thus, reading the results in each pattern from right to left and writing them in the usual manner, the four possible results are `1 0', `0 1', `0 1' and `0 0', in order from the top row to the bottom.

\begin{figure}[!bhpt]
\begin{center}
\begin{tabular}{l}
A 1 1 S 0 C 1 \\
A 1 0 S 1 C 0 \\
A 0 1 S 1 C 0 \\
A 0 0 S 0 C 0 \\
\end{tabular}
\end{center}
\caption{\small Four patterns representing a table to define the function for the addition of two one-bit numbers in binary arithmetic, with provision for the carrying out of one bit. The 
meanings of the letters are described in the text.}
\label{half_adder_definition}
\normalsize
\end{figure}

Apart from their function in conveying the meanings of the bits to human readers, 
service symbols (like `A', `S' and `C') are needed to constrain the search for possible 
alignments and, in effect, to tell the system which digits represent the `input' digits and 
which are `output'.\footnote{A distinction between `content' symbols and `code' symbols provides additional constraint but the details are not relevant to the subject of this article and need not be spelled out here.}

Figure \ref{one_bit_adder_alignment} shows the best alignment that has been found by SP61 with the pattern 
`A 0 1 S C' in New and the patterns shown in Figure \ref{half_adder_definition} in Old. (Notice that the first and last digit in each row is a row number, not part of the pattern for that row.) Generally speaking, in ICMAUS the `result' of an alignment is the symbols in patterns from Old in the alignment that have not been aligned with any symbols from New. The result in this case is `1' for the sum bit and `0' for the carry-out bit. In a very straightforward way, the matching process has functioned as a process of `table lookup' and has `selected' the pattern in Figure \ref{half_adder_definition} that yields the `correct' addition of `0' and `1' in binary arithmetic.

\begin{figure}[!bhpt]
\begin{center}
\begin{BVerbatim}
0 A 0 1 S   C   0
  | | | |   |  
1 A 0 1 S 1 C 0 1
\end{BVerbatim}
\end{center}
\caption{\small The best alignment found by SP61 with `A 0 1 S C' in New and the patterns shown in Figure \ref{half_adder_definition} in Old. The first and last digit in each row is a row number, not part of the pattern for that row.}
\label{one_bit_adder_alignment}
\normalsize
\end{figure}

In general, any function that can be represented as a look-up table can be 
represented in ICMAUS as a set of patterns and processed in essentially the same way as 
the example of one-bit addition shown here.

\subsection{COMPOSITE FUNCTIONS AND THEIR EXECUTION}

We saw in Section \ref{maths_functions} that, although all functions in ML are defined abstractly as a table of inputs and outputs, many functions do not look much like a table and are very much more compact than the table that they represent. In later parts of Section \ref{ML_structures} we have seen that the succinctness of many functions, and the succinctness of ML in general, can often be attributed to the exploitation of our three simple compression techniques: chunking-with-codes, schema-plus-correction and run-length coding.

So far, all of this seems to be reasonably straightforward. A complication, however, 
is that when two or more functions are part of a higher-level function, or when there is a 
repetition of a function two or more times, it is often necessary for the output of one 
function to become the input of another function or for the output of a function to become 
the input for another invocation of the same function. Does this passing of information 
between functions represent any kind of IC and, in particular, can it be modelled with 
ICMAUS?

The next subsection shows how, in a simple case, it can be modelled in ICMAUS. 
And the following subsection considers this issue in more general terms.

\subsubsection{\em Adding Two-Bit Numbers or Larger}\label{adding_two_bit_numbers}

Conceptually, the simplest way to add together a pair of numbers containing two or 
more bits is to apply the one-bit adder recursively to successive bits in each of the two 
numbers which are to be added. And this means that it should be possible for any carry-out 
bit from one application of the adder to become the {\it carry-in} bit for the next application of the adder in the next higher column of the addition.

In the jargon of computer engineering, a one-bit adder like the one shown in Figure \ref{half_adder_definition} is called a {\it half adder} \citep[see][pp. 83--89]{richards_1955} because it lacks any means of receiving a carry-in bit from the column immediately below. This deficiency is made up in a {\it full adder} \citep[see][pp. 89--95]{richards_1955}, a version of which is represented by the set of patterns shown in Figure \ref{full_adder_definition}. Here, the `C' bit at the beginning of each pattern is the carry-in bit and all the other bits are the same as in the half adder (Figure \ref{half_adder_definition}).

Notice that the pair of symbols `C ... \#C' marking the carry-in bit at the beginning 
of each pattern is the same as the symbols `C ... \#C' that mark the carry-out bit at the end 
of the pattern. Matching of these pairs of symbols allows the carry-out bit at any one level to become the carry-in bit at the next higher level so that values can, in effect, be passed from level to level. This can be seen in the example described next.

Consider, for example, how we may add together two three-bit numbers like `1 1 0' 
and `0 1 1'. The two numbers need to be arranged so that the bits for each level lie together 
and are marked in the same way as in the New pattern in Figure \ref{one_bit_adder_alignment}. The order needs to be reversed to take account of the way the SP61 model searches for alignments.\footnote{In principle, SP61 can search for good alignments by processing all the symbols in New together. But in practice, some types of good alignment (including the ones described here) are much easier to find if the symbols from New are processed one at a time in left-to-right order. And, in addition, it is natural to process columns in sequence from the lowest order to the highest.} Thus the sequence of symbols for these two numbers should be: `A 1 0 S A 1 1 S A 0 1 S'. However, the symbols `C 0 \#C' need to be added at the beginning of New to ensure that the carry-in bit for the first addition is `0'. So the final result is `C 0 \#C A 1 0 S A 1 1 S A 0 1 S'.

The best alignment found by SP61 with `C 0 \#C A 1 0 S A 1 1 S A 0 1 S' in New 
and the patterns shown in Figure \ref{full_adder_definition} in Old is shown in Figure \ref{three_bit_addition_alignment}. Ignoring the `service' symbols, the sum of the two numbers can be read from this alignment as the three `sum' bits (each one between an `S' and a `C' column), with the order reversed, and preceded on the left by the last carry-out bit on the right of the alignment. The overall result is `1 0 0 1'. Readers will verify easily enough that this is the correct sum of the binary numbers `1 1 0' and `0 1 1'.

\begin{figure}[!bhpt]
\begin{center}
\begin{tabular}{l}
C 1 \#C A 1 1 S 1 C 1 \#C \\
C 1 \#C A 0 1 S 0 C 1 \#C \\
C 1 \#C A 1 0 S 0 C 1 \#C \\
C 0 \#C A 1 1 S 0 C 1 \#C \\
C 1 \#C A 0 0 S 1 C 0 \#C \\
C 0 \#C A 0 1 S 1 C 0 \#C \\
C 0 \#C A 1 0 S 1 C 0 \#C \\
C 0 \#C A 0 0 S 0 C 0 \#C \\
\end{tabular}
\end{center}
\caption{\small Eight patterns representing a table to define the function for the addition of two one-bit numbers in binary arithmetic, with provision for the carrying in of a digit as well as the carrying out of a digit. The meanings of the symbols are described in the text.}
\label{full_adder_definition}
\normalsize
\end{figure}

\begin{figure}[!bhpt]
\fontsize{10.00pt}{12.00pt}
\begin{center}
\begin{BVerbatim}
0 C 0 #C A 1 0 S          A 1 1 S          A 0 1 S          0
  | | |  | | | |          | | | |          | | | |         
1 C 0 #C A 1 0 S 1 C 0 #C | | | |          | | | |          1
                   | | |  | | | |          | | | |         
2                  C 0 #C A 1 1 S 0 C 1 #C | | | |          2
                                    | | |  | | | |         
3                                   C 1 #C A 0 1 S 0 C 1 #C 3
\end{BVerbatim}
\end{center}
\caption{\small The best alignment found by SP61 with the pattern `C 0 \#C A 1 0 S A 1 1 S A 0 1 S' in New and the patterns shown in Figure \ref{full_adder_definition} in Old.}
\label{three_bit_addition_alignment}
\normalsize
\end{figure}

\subsubsection{\em Other Composite Functions}

Compared with conventional systems, the method just shown for passing 
information from one invocation of a function to another is somewhat inflexible. This is 
because it requires the `output' variable to have the same name as the `input' variable. In 
conventional systems there is a variety of notations for showing that the output of one 
function becomes the input to another and none of them are restricted in this way. 

Modelling these more flexible methods in ICMAUS probably requires an ability to 
model `scoping' of variables and, possibly, an ability for `learning'. Both these things are beyond what the SP61 model can do as it has been developed to date and discussion of these issues would take up too much space in this article.

\subsection{PROPOSITIONAL LOGIC AND IC}

Figure \ref{xor_definition} shows four patterns representing a `truth table' for the `exclusive OR' 
logical proposition ($p \vee q$. ``Either $p$ is true or $q$ is true but they are not both true''). Here, `1' means TRUE and `0' means FALSE, the first two digits in each pattern are the truth values of $p$ and $q$ (individually) and, for each pattern, the last digit, between the symbols `R .. \#R' , is the truth value of $p \vee q$.

Readers will notice that the XOR truth table is very similar to the definition of 
one-bit addition shown in Figure \ref{half_adder_definition}, not merely in the arrangement of the digits but, more importantly, in the fact that both things can be represented as a simple lookup table that can be represented by a set of patterns. It should be no surprise that the XOR function can be evaluated by SP61 in essentially the same way as with the function for one-bit addition.

\begin{figure}[!bhpt]
\begin{center}
\begin{tabular}{l}
XOR 1 1 R 0 \#R \\
XOR 0 1 R 1 \#R \\
XOR 1 0 R 1 \#R \\
XOR 0 0 R 0 \#R \\
\end{tabular}
\end{center}
\caption{\small Four patterns representing the `exclusive OR' (XOR) truth table as described in the text.}
\label{xor_definition}
\normalsize
\end{figure}

Again, one might expect to be able to evaluate combinations of simple propositions 
in ICMAUS in essentially the same way as was done for the addition of numbers containing 
two or more bits (Section \ref{adding_two_bit_numbers}). As before, this can be achieved with the SP61 model in its current stage of development if outputs are made to match inputs. This can be seen in the definition of `NOTXOR' shown in Figure \ref{notxor_definition}. Here, the alternative `outputs' of the XOR proposition are contained in the field `A ... \#A' and the same two symbols are used to mark the possible `inputs' of the truth table for NOT (that changes `1' to `0' and {\it vice versa}).

\begin{figure}[!bhpt]
\begin{center}
\begin{tabular}{l}
NOTXOR XOR NOT R \#R \#NX \\
XOR 1 1 A 0 \#A \\
XOR 1 0 A 1 \#A \\
XOR 0 1 A 1 \#A \\
XOR 0 0 A 0 \#A \\
A 1 \#A NOT R 0 \#R \\
A 0 \#A NOT R 1 \#R \\
\end{tabular}
\end{center}
\caption{\small A set of patterns defining the `not XOR' proposition as described in the text.}
\label{notxor_definition}
\normalsize
\end{figure}

\sloppy{Figure \ref{notxor_alignment} shows the best alignment found by SP61 with the pattern `NOTXOR 0 1 A \#NX' in New and the patterns from Figure \ref{notxor_definition} in Old. In the alignment, the digits `0' and `1' in New yield `A 1 \#A' as the `result' of the XOR component. These three symbols match the first three symbols in the first of the two `NOT' patterns and this yields `R 0 \#R as the overall `result' of the NOTXOR proposition.}

\begin{figure}[!bhpt]
\begin{center}
\begin{BVerbatim}
0 NOTXOR     0 1 A                 #NX 0
    |        | | |                  | 
1   |    XOR 0 1 A 1 #A             |  1
    |     |      | | |              | 
2   |     |      A 1 #A NOT R 0 #R  |  2
    |     |              |  |   |   | 
3 NOTXOR XOR            NOT R   #R #NX 3
\end{BVerbatim}
\end{center}
\caption{\small The best alignment found by SP61 with `NOTXOR 0 1 A \#NX' in New and the 
patterns shown in Figure \ref{notxor_definition} in Old.}
\label{notxor_alignment}
\normalsize
\end{figure}

As in the case of adding numbers containing two or more digits, the ICMAUS 
patterns that have been shown are somewhat contrived owing to the need to make the 
output at one stage match the input at the next stage. It is anticipated that more elegant 
solutions may be found when the SP61 has been further developed.

\subsection{TRUE AND FALSE IN ICMAUS}\label{TRUE_FALSE_ICMAUS}

Although readers may accept that, in broad terms, the matching and unification of patterns in ICMAUS can model the kinds of inferences that can be made in propositional logic, an objection may be raised that, in propositional logic {\it per se}, every inference has only two values, TRUE or FALSE, whereas in ICMAUS, there may be and often are alternative alignments for a given inference that suggests that there are uncertainties attaching to the inferences that the system makes.

ICMAUS applies most naturally to probabilistic reasoning \citep[see][]{wolff_1999_prob} and this suggests that `exact' reasoning may be seen as reasoning where probabilities are very close to 0 or 1. In itself, this idea is not satisfactory owing to uncertainties about the precise meaning of `close to' in this context. A better idea is to use explicit TRUE and FALSE values in `output' fields (or `0' and `1' values as in the examples we have been considering), coupled with a focus always on the `best' alignment for any set of patterns. If alignments are created and interpreted in this way then the effect of two-valued reasoning can be achieved within the ICMAUS framework without ambiguity. An example is shown in the next subsection.

\subsection{SYLLOGISTIC REASONING}\label{syllogistic_reasoning}

Consider the following text book example of a {\it modus ponens} syllogism:

\begin{enumerate}

\item All humans are mortal.

\item Socrates is human.

\item Therefore, Socrates is mortal.

\end{enumerate}

\noindent In logical notation, this may be expressed as:

\begin{enumerate}

\item $\forall x: human(x) \Rightarrow mortal(x)$.

\item $human(Socrates)$.

\item $\therefore mortal(Socrates)$.

\end{enumerate}

The traditional `strict' interpretation in logic of syllogistic reasoning requires 
explicit values (TRUE and FALSE) for propositions, in much the same style as the `output' 
fields of the XOR, NOT and NOTXOR functions shown above. 

In the ICMAUS framework, the first of the propositions just shown may be expressed with the pattern `X \#X human true $\Rightarrow$ mortal true'. In the manner of Skolemization, the variable `X \#X' in the pattern represents anything at all and may thus be seen to be universally quantified. The scope of the variable may be seen to embrace the entire pattern, without the need for it to be repeated. The symbol `$\Rightarrow$' in the ICMAUS pattern serves simply as a separator between `human true' and `mortal true'.

If this pattern is included (with other patterns) in Old and if New contains the pattern `Socrates human true $\Rightarrow$' (corresponding to `human(Socrates)' and, in effect, a request to discover what that proposition implies), SP61 finds one alignment (shown in Figure \ref{truth_value_socrates_alignment}) that encodes all the symbols in New. After unification to form the pattern `X Socrates \#X human true $\Rightarrow$ mortal true', we may read the alignment as a statement that because it is true that Socrates is human it is also true that Socrates is mortal. In this reading, we ignore the `service' symbols, `X' and `\#X'.

\begin{figure}[!bhpt]
\begin{center}
\begin{BVerbatim}
0   Socrates    human true =>             0
       |          |    |   |             
1 X    |     #X human true => mortal true 1
  |    |     |                           
2 X Socrates #X                           2
\end{BVerbatim}
\end{center}
\caption{The best alignment found by SP61 with `Socrates human true $\Rightarrow$' in New and patterns in Old that include `X \#X human true $\Rightarrow$ mortal true'.}
\label{truth_value_socrates_alignment}
\end{figure}

Notwithstanding the probabilistic nature of the ICMAUS framework, there is no ambiguity in an alignment like this because it is the only alignment formed that encodes all the symbols in New. Given that we restrict our attention to alignments of this kind, it is entirely possible for the framework to yield deductions that are either TRUE or FALSE and not something in between.

\subsection{PROBABILISTIC REASONING}

Apart from the `strict' logical interpretation of examples like the one above, there 
is a more relaxed `everyday' interpretation that allows statements and inferences to have 
levels of confidence or `probabilities'. This kind of interpretation would view a proposition 
like ``All humans are mortal'' to be an inductive generalisation about the world which, like 
all inductive generalisations, might turn out to be wrong. Likewise, we may not be totally 
confident that Socrates is human. As a consequence of uncertainties like these, the 
proposition that ``Socrates is mortal'' may have some doubt attaching to it.

A relatively full discussion of the way the ICMAUS framework may be applied to probabilistic reasoning (including such things as default reasoning, nonmonotonic reasoning and `explaining away') may be found in \citet{wolff_1999_prob}.

\subsection{CHAINS OF REASONING}

By now it should be clear that more complex kinds of reasoning, including `chains' 
of reasoning, may be built up from simpler kinds, by arranging patterns so that, by matching 
and unification of symbols, the `output' of one stage can determine the `input' of the next 
stage. For example, simple `if ... then' rules of the kind that are the stock in trade of expert systems may be modelled in ICMAUS with patterns that associate {\it conditions} with 
{\it actions} or {\it consequences} (e.g., ``If smoke, then fire'', ``If litmus paper is red, then solution is acid'').

In the simplest schematic terms, a rule like `IF A THEN B' may, in ICMAUS, be 
simplified to the association of `A' and `B' in a pattern like `A B'. Then an argument like 
``A means C because A means B and B means C'' can be modelled by the alignment of an 
initial `A' to the first `A' in `A B', followed by the alignment of the two instances of `B' in the patterns `A B' and `B C'.

Readers wishing to pursue these ideas will find other examples and relevant 
discussion in \citet{wolff_1999_prob}. For the purposes of this article, the key point here is the way in which inferences may be understood in terms of IC arising from the alignment of patterns and parts of patterns.

\section{\bf Discussion}\label{discussion_conclusion}

Before concluding this article, let us consider briefly a few issues that 
relate to the ideas that have been presented. There is insufficient space to discuss these issues fully---this section merely flags their relevance and indicates some lines of thinking.

\subsection{THE PARADOX OF `DECOMPRESSION BY COMPRESSION'}\label{paradox}

Given the thesis that ML may be understood in terms of IC, it is reasonable to ask 
how it is possible, as it clearly is, for decompression to be achieved by mathematical 
processes of calculation and inference. Newton's formula, $s = (gt^{2}) / 2$, may be seen as a 
compressed representation of the table relating distances travelled by an object falling 
freely from rest in a given time. But, with ordinary mathematical processes of calculation, 
the formula may be used to recreate large numbers of rows in that table and thus 
decompress the information.

\sloppy{This paradox and how it may be resolved is discussed fairly fully in \citet{wolff_2000}, with examples.} In essence, the idea is to replace the original New with an encoded, compressed version and run the system again. Generally, it will recreate the original alignment (except for the substitution in New). Unification of the new alignment yields the elements of the original New, in the same order, together with `service' symbols that can be ignored.

The reason that there is no real paradox is that the `query' pattern used to retrieve 
information is new information. It is the compression of this new information by the 
unification of patterns that match each other that achieves the retrieval of the previously
compressed information. This retrieval of information achieves the effect of 
decompression.

\subsection{COULD IC BE A FOUNDATION FOR ML?}

Rather than suppose that mathematics is founded on logic or {\it vice versa}, this article 
has suggested that mathematics, logic and related disciplines may be seen to be founded on IC, itself the product of primitive operations of matching and unification of patterns.

Possible objections to this thesis are that:

\begin{itemize}

\item Concepts of information that are the basis of methods for IC are couched in 
mathematical terms and

\item Mathematics is used in methods for IC, including the ICMAUS model.

\end{itemize}

A possible answer to the first point is that the Hartley-Shannon and AIT concepts 
of information are themselves founded on a more primitive idea of discriminable changes (patterns of variation) across one or more dimensions and that the role of mathematics is to measure information, not define the concept itself. Mathematics is useful as a way of giving precision to the measurement of information but it is not the only way in which information may be measured (see below).

In a similar vein, an answer to the second point is that mathematics can be useful in 
methods for IC but that the core operations of matching and unifying patterns do not, in themselves, use mathematical constructs.

Rather less certainly, it may be said that, although a concept like `frequency' (which 
has an important role in discriminating amongst alternative ways in which patterns may be 
matched and unified) may conveniently be expressed in the form of numbers (that are 
clearly mathematical) the concept may equally well be expressed as a non-numeric 
`strength' of some kind of chemical or physical variable in a brain or computer. Much 
depends, of course, on what does or does not count as being `mathematical'.

In general, there are alternatives to mathematics in the measurement of information, 
in compression of information and in the measurement of frequencies of occurrence. The 
proof of these assertions is that all of these things can be and often are done by brains and 
nervous systems in humans and other animals without the assistance of mathematical knowledge. It would be perverse to say that a fox's ability to learn where rabbits are most likely to be found is a mathematical ability. Capabilities of this sort may be less accurate than those that exploit mathematical techniques but they are useful for everyday living and probably essential for survival of species through millions of years.

\subsection{ML AND HUMAN PSYCHOLOGY}

Section \ref{brains_and_ns} briefly reviewed some of the writings and research relating to the idea that much of learning, perception and cognition in humans (and animals) may be seen to be founded on information compression. In that section and later it was suggested that 
concepts in ML may be seen to be an outgrowth or development of this fundamental feature 
of human cognition.

Here we simply note this suggestion as a matter for further consideration in the 
future. Relevant issues include the nature of the evidence that one might recognise as being 
relevant to the acceptance or rejection of the hypothesis, whether or not the hypothesis 
yields insights into ML that might otherwise be missed and, in a similar way, whether the 
hypothesis leads to new insights into human learning, perception and cognition.

\subsection{MATHEMATICS AND ITS APPLICABILITY}\label{maths_science}

The arguments and examples in this article may provide some solution to the 
`mystery' of why mathematics is such a good ``language of science''. As Barrow says 
\citeyearpar[p. 247]{barrow_1992}, ``Science is, at root, just the search for compression in the world''. If the arguments in this article are accepted, then mathematics is a good medium for the expression of laws and other regularities in the world because mathematics provides the kinds of mechanisms that are standardly used for IC, with matching and unification of patterns at their core.

Notice that the usefulness of IC in this connection is not merely in reducing the volume of scientific information (which may facilitate storage or transmission of that information). At least as important is the intimate connection that appears to exist between IC and the inductive prediction of the future from the past (Section \ref{inf_MLE_PR}, above, and Section \ref{inductive_prediction}, below).

\subsubsection{\em The World is Compressible}\label{the_world_is_compressible}

Barrow goes on to say [{\it ibid.}, p. 247], ``... the world is surprisingly compressible and the success of mathematics in describing its workings is a manifestation of that 
compressibility.'' 

It is interesting to reflect on what things would be like if the world (and 
the universe) were {\it not} compressible. In accordance with AIT principles, such a world 
would be totally random. This would mean that there would be no `structure' of any kind, 
including the world itself, people, brains, computers and, indeed, the structures of 
mathematics, logic and related disciplines. So all of the things we have been discussing, 
including mathematics and science, would not exist. Science and the expression of 
scientific truths in mathematics can only exist {\it because} the world is compressible!

It is true that, out of all the possible patterns of information, compressible patterns 
(containing regularities) are relatively rare \citep{li_vitanyi_1997}. If surprisingness depends on rarity, then it is surprising that the world is compressible. But since, in accordance with the anthropic principle, we can only exist in a compressible world that accommodates persistent structures like ourselves, it is {\it not} surprising that the world is compressible!

\subsubsection{\em The Scope of Mathematics in Science}

Mathematics is useful in many areas of science but it has limitations. The theory of 
evolution by natural selection, which is by any standards a SUPERB theory (to use Penrose's word and emphasis) is not, fundamentally, mathematical. Mathematics can certainly be used to 
analyse aspects of evolution but the theory in itself does not use any mathematical constructs (unless one adopts a definition of mathematics that is so wide that it becomes meaningless). Again, there are many aspects of the world---like distributions of plants and animals or patterns of political beliefs---that seem to be fundamentally `messy' and not amenable to `neat' compression into mathematical formulae.

In short, mathematics is a good language of science but it is not the only language 
for expressing regularities in the world. And there are many aspects of the world that 
have not yet been compressed into compact formulae or any other succinct representation and seem unlikely to yield to that treatment in the future.

\subsection{ML AND INDUCTIVE PREDICTION OF THE FUTURE 
FROM THE PAST}\label{inductive_prediction}

If the main thesis of this article is accepted, then it is accepted that ML is founded 
on the same principles that appear to underlie concepts of probability and probabilistic 
inductive inference. As we saw in Section \ref{inf_MLE_PR} and later: there is a very close connection between IC and concepts of probability; the ICMAUS framework is based on MLE 
principles; and \citet{solomonoff_1986} has argued that the great majority of problems in science {\it and} mathematics can be seen as either `machine inversion' problems or `time limited optimization' problems, and that both kinds of problems can be solved by inductive inference using the MLE principle.

Of course, there are parts of mathematics to do with probability and there are types 
of logic that are probabilistic. But, for many mathematicians and logicians, the idea that 
mathematics or logic are founded on the same principles that underlie concepts of 
probability and probabilistic inductive inference may be hard to accept. It is, after all, the 
highly predictable `clockwork' nature of these disciplines that gives them their appeal.

There is no space here to discuss these issues fully. Five points are noted here briefly:

\begin{itemize}

\item Non-determinism and associated uncertainties can appear in the Turing model 
itself, depending on the transition function built into the model. Likewise, the Post 
Canonical System model can exhibit non-determinism, depending on the 
productions provided in the model.

\item Section \ref{TRUE_FALSE_ICMAUS} described briefly how the effect of two-valued logic might be obtained in a system that is otherwise geared to probabilistic reasoning.

\item As we saw in Section \ref{syllogistic_reasoning}, some kinds of syllogistic argument can be seen to have inductive origins.

\item That mathematics may not, fundamentally, be quite as `clockwork' as is commonly 
supposed has been suggested by \citet[p. 80]{chaitin_1988}: ``{\it I have recently been able 
to take a further step along the path laid out by G\"{o}del and Turing. By translating a 
particular computer program into an algebraic equation of a type that was familiar 
even to the ancient Greeks, $I$ have shown that there is randomness in the branch of 
pure mathematics known as number theory. My work indicates that---to borrow 
Einstein's metaphor---God sometimes plays dice with whole numbers.}''

\item The idea that we can only exist in a compressible world (Section \ref{the_world_is_compressible}, above) seems to provide an answer to the problem of finding a rational basis for inductive reasoning. Of course, we cannot justify such reasoning by saying that it has always worked in the past---because this simply invokes the principle we are trying to justify. A better answer seems to be as follows:

\begin{itemize}

\item The assumption that the future will repeat the past is equivalent to the assumption that there is information redundancy between past and future.

\item If we make that assumption and it turns out to be true, then we will reap the benefits of inductive reasoning---remembering where food or shelter can be found, avoiding dangers, and so on.

\item If we make that assumption and it turns out to be false, the universe will have exploded into total randomness (or exploded and reformed into {\it totally} unrecognisable patterns) and we are dead anyway!

\item In short, there is no downside risk to the assumption that the past is a guide to the future and a very substantial upside if it turns out to be true. Hence, it is a rational assumption to make.

\noindent This line of thinking is similar to Reichenbach's ``pragmatic vindication'' of induction \citep[see, for example,][]{salmon_1991}.

\end{itemize}

\end{itemize}

\subsection{RECONCILING THE INEVITABILITY OF MATHEMATICS WITH ITS APPLICABILITY}

\citet{potter_2000} writes that the chief philosophical problem posed by arithmetic is ``that of reconciling its inevitability---the impression we have that it is true no matter what---with its applicability.'' (p. 1). Presumably, the same can be said about mathematics as a whole and, indeed, logic.

This section has already suggested a reason why ML is so applicable: because it is an embodiment of standard techniques for information compression and, as such, it provides a convenient means of describing the world in a succinct manner and, perhaps more importantly, it gives us a handle on the inductive prediction of the future from the past.

What about the `inevitability' of arithmetic and other aspects of ML? The first point in this connection is that, as was noted above, the work of G\"{o}del, Turing, Chaitin and others has shown that arithmetic, at least, has a probabilistic aspect (``God sometimes plays dice with whole numbers''): it may not be as true as is commonly assumed that 2 + 2 = 4, no matter what. That said, any `fuzziness' in arithmetic is not apparent in its everyday applications and it is still pertinent to ask why the laws of arithmetic are amongst the most certain things that we can know.

The view of ML that has been sketched in this article suggests possible reasons for the certainties of arithmetic and, indeed, the uncertainties too:

\begin{itemize}

\item With the simple kinds of patterns that seem to form the basis of natural numbers (Section \ref{generation_recognition_numbers}, above), ICMAUS often yields an alignment that stands out very clearly as better in terms of compression than any of the alternative alignments---because simple patterns impose relatively tight constraints on the number of alternative ways in which symbols may be matched and unified. This phenomenon may be responsible for our subjective impression of the inevitability of arithmetic truths. 

\item The uncertainties and paradoxes of ML (G\"{o}del's theorems, the Halting Problem, Chaitin's work, Russell's paradox, etc) seem all to arise from the infinite regress of recursion that can occur in ML systems. This sits comfortably with the ICMAUS view of ML because recursion can easily arise in the matching and unification of patterns (as we saw in Section \ref{generation_recognition_numbers}).

\end{itemize}

\section{\bf Conclusion}

I hope that the view of mathematics, logic and related disciplines that is presented in this article will prove interesting to readers. It is a long way from being any kind of `proof' that ML is based on IC but may provide a perspective that stimulates new thinking about the nature of mathematics, logic and related disciplines, and new insights into those disciplines.

\section*{Acknowledgements}

I am very grateful to Peter Apostoli, Michele Friend, Alex Paseau, Oliver Schulte and Robert Thomas for constructive comments on earlier drafts of this paper. I am also grateful for discussion of some of the ideas in this article with Tim Porter and Chris Wensley.


\begin{thebibliography}{46}
\expandafter\ifx\csname natexlab\endcsname\relax\def\natexlab#1{#1}\fi

\bibitem[Attneave(1954)]{attneave_1954}
F.~Attneave.
\newblock Some informational aspects of visual perception.
\newblock {\em Psychological Review}, 61:\penalty0 183--193, 1954.

\bibitem[Barlow(1959)]{barlow_1959}
H.~B. Barlow.
\newblock Sensory mechanisms, the reduction of redundancy, and intelligence.
\newblock In {HMSO}, editor, {\em The Mechanisation of Thought Processes},
  pages 535--559. Her Majesty's Stationery Office, London, 1959.

\bibitem[Barlow(1969)]{barlow_1969}
H.~B. Barlow.
\newblock Trigger features, adaptation and economy of impulses.
\newblock In K.~N. Leibovic, editor, {\em Information Processes in the Nervous
  System}, pages 209--230. Springer, New York, 1969.

\bibitem[Barlow(1997)]{barlow_1997}
H.~B. Barlow.
\newblock The knowledge used in vision and where it comes from.
\newblock {\em Philosophical Transactions of the Royal Society London B},
  352:\penalty0 1141--1147, 1997.

\bibitem[Barlow(2001)]{barlow_2001_bbs}
H.~B. Barlow.
\newblock The exploitation of regularities in the environment by the brain.
\newblock {\em Behavioural and Brain Sciences}, 24\penalty0 (4):\penalty0
  602--607 and 748--749, 2001.

\bibitem[Barrow(1992)]{barrow_1992}
J.~D. Barrow.
\newblock {\em Pi in the Sky}.
\newblock Penguin Books, Harmondsworth, 1992.

\bibitem[Birtwistle et~al.(1973)Birtwistle, Dahl, Myhrhaug, and
  Nygaard]{birtwistle_1973}
G.~M. Birtwistle, O-J Dahl, B.~Myhrhaug, and K.~Nygaard.
\newblock {\em Simula Begin}.
\newblock Studentlitteratur, Lund, 1973.

\bibitem[Boolos and Jeffrey(1980)]{boolos_jeffrey_1980}
G.~Boolos and R.~Jeffrey.
\newblock {\em Computability and Logic}.
\newblock Cambridge University Press, Cambridge, 1980.

\bibitem[Chaitin(1988)]{chaitin_1988}
G.~J. Chaitin.
\newblock Randomness in arithmetic.
\newblock {\em Scientific American}, 259\penalty0 (1):\penalty0 80--85, 1988.

\bibitem[Chater(1996)]{chater_1996}
N.~Chater.
\newblock Reconciling simplicity and likelihood principles in perceptual
  organisation.
\newblock {\em Psychological Review}, 103\penalty0 (3):\penalty0 566--581,
  1996.

\bibitem[Chater(1999)]{chater_1999}
N.~Chater.
\newblock The search for simplicity: a fundamental cognitive principle?
\newblock {\em Quarterly Journal of Experimental Psychology}, 52 A\penalty0
  (2):\penalty0 273--302, 1999.

\bibitem[Chomsky(1957)]{chomsky_1957}
N.~Chomsky.
\newblock {\em Syntactic Structures}.
\newblock Mouton, The Hague, 1957.

\bibitem[Cover and Thomas(1991)]{cover_thomas_1991}
T.~M. Cover and J.~A. Thomas.
\newblock {\em Elements of Information Theory}.
\newblock John Wiley, New York, 1991.

\bibitem[Craig(1998)]{craig_1998}
E.~Craig, editor.
\newblock {\em Routledge Encyclopedia of Philosophy}, volume 1--10.
\newblock Routledge, London, 1998.

\bibitem[Devlin(1991)]{devlin_1991}
K.~Devlin.
\newblock {\em Logic and Information}.
\newblock Cambridge University Press, Cambridge, 1991.

\bibitem[Devlin(1997)]{devlin_1997}
K.~Devlin.
\newblock {\em Mathematics: The Science of Patterns}.
\newblock Scientific American Library, New York, 1997.

\bibitem[Epstein and Carnielli(1989)]{epstein_carnielli_1989}
R.~L. Epstein and W.~A. Carnielli.
\newblock {\em Computability: Computable Functions, Logic, and the Foundations
  of Mathematics}.
\newblock Wadsworth \& Brooks, Pacific Grove, 1989.

\bibitem[Eves(1990)]{eves_1990}
H.~Eves.
\newblock {\em Foundations and Fundamental Concepts of Mathematics}.
\newblock Dover Publications, New York, 1990.

\bibitem[Garner(1974)]{garner_1974}
W.~R. Garner, editor.
\newblock {\em The Processing of Information and Structure}.
\newblock Lawrence Erlbaum, Hillsdale, NJ, 1974.

\bibitem[Hart(1996)]{hart_1996}
W.~D. Hart.
\newblock {\em The Philosophy of Mathematics}.
\newblock Oxford University Press, New York, 1996.

\bibitem[Hersh(1997)]{hersh_1997}
R.~Hersh.
\newblock {\em What is Mathematics Really?}
\newblock Vintage, London, 1997.

\bibitem[Li and Vit\'{a}nyi(1997)]{li_vitanyi_1997}
M.~Li and P.~Vit\'{a}nyi.
\newblock {\em An Introduction to Kolmogorov Complexity and Its Applications}.
\newblock Springer-Verlag, New York, 1997.

\bibitem[Mandelbrot(1993)]{plaisted_1993}
B.~Mandelbrot.
\newblock Linguistique statistique macroscopique.
\newblock In D.~M. Gabbay, C.~J. Hogger, and J.~A. Robinson, editors, {\em
  Handbook of Logic in Artificial Intelligence and Logic Programming: Logical
  Foundations}, volume~1, pages 1--78. Oxford University Press, Oxford, 1993.

\bibitem[Oldfield(1954)]{oldfield_1954}
R.~C. Oldfield.
\newblock Memory mechanisms and the theory of schemata.
\newblock {\em British Journal of Psychology}, 45:\penalty0 14--23, 1954.

\bibitem[Penrose(1989)]{penrose_1989}
R.~Penrose.
\newblock {\em The Emperor's New Mind}.
\newblock Oxford University Press, Oxford, 1989.

\bibitem[Post(1943)]{post_1943}
E.~L. Post.
\newblock Formal reductions of the general combinatorial decision problem.
\newblock {\em American Journal of Mathematics}, 65:\penalty0 197--268, 1943.

\bibitem[Potter(2000)]{potter_2000}
M.~Potter.
\newblock {\em Reason's Nearest Kin: Philosophies of Arithmetic from Kant to
  Carnap}.
\newblock Oxford University Press, Oxford, 2000.

\bibitem[Resnik(1997)]{resnik_1997}
M.~D. Resnik.
\newblock {\em Mathematics as a Science of Patterns}.
\newblock Clarendon Press, Oxford, 1997.

\bibitem[Richards(1955)]{richards_1955}
R.~K. Richards.
\newblock {\em Arithmetic Operations in Digital Computers}.
\newblock Van Nostrand, Princeton, NJ, 1955.

\bibitem[Salmon(1991)]{salmon_1991}
W.~Salmon.
\newblock Hans reichenbach's vindication of induction.
\newblock {\em Erkenntnis}, 35:\penalty0 99--122, 1991.

\bibitem[Sankoff and Kruskall(1983)]{sankoff_kruskall_1983}
D.~Sankoff and J.~B. Kruskall.
\newblock {\em Time Warps, String Edits, and Macromolecules: the Theory and
  Practice of Sequence Comparisons}.
\newblock Addison-Wesley, Reading, MA, 1983.

\bibitem[Shannon and Weaver(1949)]{shannon_weaver_1949}
C.~E. Shannon and W.~Weaver.
\newblock {\em The Mathematical Theory of Communication}.
\newblock University of Illinois Press, Urbana, 1949.

\bibitem[Shapiro(2000)]{shapiro_2000}
S.~Shapiro.
\newblock {\em Philosophy of Mathematics}.
\newblock Oxford University Press, New York, 2000.

\bibitem[Solomonoff(1986)]{solomonoff_1986}
R.~J. Solomonoff.
\newblock The application of algorithmic probability to problems in artificial
  intelligence.
\newblock In L.~N. Kanal and J.~F. Lemmer, editors, {\em Uncertainty in
  Artificial Intelligence}, pages 473--491. Elsevier Science, North-Holland,
  1986.

\bibitem[Storer(1988)]{storer_1988}
J.~A. Storer.
\newblock {\em Data Compression: Methods and Theory}.
\newblock Computer Science Press, Rockville, Maryland, 1988.

\bibitem[Sudkamp(1988)]{sudkamp_1988}
T.~A. Sudkamp.
\newblock {\em Languages and Machines}.
\newblock Addison-Wesley, Reading, Mass., 1988.

\bibitem[{von B{\'e}k{\'e}sy}(1967)]{von_bekesy_1967}
G.~{von B{\'e}k{\'e}sy}.
\newblock {\em Sensory Inhibition}.
\newblock Princeton University Press, Princeton, NJ, 1967.

\bibitem[Watanabe(1972)]{watanabe_article_1972}
S.~Watanabe.
\newblock Pattern recognition as information compression.
\newblock In S.~Watanabe, editor, {\em Frontiers of Pattern Recognition}.
  Academic Press, New York, 1972.

\bibitem[Wolff(1982)]{wolff_1982}
J.~G. Wolff.
\newblock Language acquisition, data compression and generalization.
\newblock {\em Language \& Communication}, 2:\penalty0 57--89, 1982.

\bibitem[Wolff(1988)]{wolff_1988}
J.~G. Wolff.
\newblock Learning syntax and meanings through optimization and distributional
  analysis.
\newblock In Y.~Levy, I.~M. Schlesinger, and M.~D.~S. Braine, editors, {\em
  Categories and Processes in Language Acquisition}, pages 179--215. Lawrence
  Erlbaum, Hillsdale, NJ, 1988.

\bibitem[Wolff(1993)]{wolff_1993}
J.~G. Wolff.
\newblock Computing, cognition and information compression.
\newblock {\em AI Communications}, 6\penalty0 (2):\penalty0 107--127, 1993.

\bibitem[Wolff(1999{\natexlab{a}})]{wolff_1999_comp}
J.~G. Wolff.
\newblock {`Computing'} as information compression by multiple alignment,
  unification and search.
\newblock {\em Journal of Universal Computer Science}, 5\penalty0
  (11):\penalty0 777--815, 1999{\natexlab{a}}.

\bibitem[Wolff(1999{\natexlab{b}})]{wolff_1999_prob}
J.~G. Wolff.
\newblock Probabilistic reasoning as information compression by multiple
  alignment, unification and search: an introduction and overview.
\newblock {\em Journal of Universal Computer Science}, 5\penalty0 (7):\penalty0
  418--462, 1999{\natexlab{b}}.

\bibitem[Wolff(2000)]{wolff_2000}
J.~G. Wolff.
\newblock Syntax, parsing and production of natural language in a framework of
  information compression by multiple alignment, unification and search.
\newblock {\em Journal of Universal Computer Science}, 6\penalty0 (8):\penalty0
  781--829, 2000.

\bibitem[Wolff(2001)]{wolff_2001_igpl}
J.~G. Wolff.
\newblock Information compression by multiple alignment, unification and search
  as a framework for human-like reasoning.
\newblock {\em Logic Journal of the IGPL}, 9\penalty0 (1):\penalty0 205--222,
  2001.
\newblock First published in the {\it Proceedings of the International
  Conference on Formal and Applied Practical Reasoning (FAPR 2000)}, September
  2000, ISSN 1469--4166.

\bibitem[Zipf(1949)]{zipf_1949}
G.~K. Zipf.
\newblock {\em Human Behaviour and the Principle of Least Effort}.
\newblock Hafner, New York, 1949.

\end{thebibliography}
\end{document}